\documentclass[12pt]{article}
\usepackage{dbw,amsmath,amssymb,amsfonts,epsfig,xspace,pstricks}
\usepackage[round]{natbib}
\pstheader{error.pro}
\pstheader{graph.pro}
\pstheader{domino.pro}
\pstheader{house.pro}

\newenvironment{remark}{\noindent{\bf Remark}\hspace*{1em}}{\bigskip}

\newcommand{\sref}[1]{\S~\ref{sec:#1}}
\newcommand{\tref}[1]{Theorem~\ref{thm:#1}}

\newcommand{\fref}[1]{Figure~\ref{fig:#1}}

\newcommand{\D}{\mathbb{D}}

\def\Cox{\hfill \square\medskip}
\newcommand{\cH}{{\cal H}}  
\renewcommand{\H}{{\cal H}}   
\newcommand{\R}{{\mathbb R}}

\begin{document}
\pagestyle{myheadings}
\markright{\sc the electronic journal of combinatorics 7 (2000), \#R00\hfill}
\thispagestyle{empty}

\title{\vspace*{-50pt}
Trees and Matchings}
\author{
\begin{tabular}{c}
Richard W. Kenyon \\
 \small Laboratoire de Topologie \\ \small Universit\'e Paris-Sud\\
{\small\texttt{kenyon\char64topo.math.u-psud.fr}}
\end{tabular}
\hspace{0.68em}
\begin{tabular}{c}
James G. Propp\thanks{Supported by NSA grant MDA904-92-H-3060,
NSF grant DMS 92-06374, and a grant from the MIT class of 1922.} \\
 \small University of Wisconsin \\
 \small Madison, Wisconsin\\
 \small\texttt{propp\char64math.wisc.edu}
\end{tabular}
\hspace{0.68em}
\begin{tabular}{c}
David B. Wilson \\
 \small Microsoft Research\\
 \small Redmond, Washington\\
 \small\texttt{dbwilson\char64alum.mit.edu}
\end{tabular}
}
\date{\small Submitted March 3, 1999; Accepted July 3, 1999; Revised March 30, 2000}
\maketitle

\vspace*{-20pt}
\begin{abstract}
In this article, Temperley's bijection between spanning trees of
the square grid on the one hand, and perfect matchings 
(also known as dimer coverings) of the
square grid on the other, is extended to the setting of general 
planar directed (and undirected) graphs, where edges carry nonnegative 
weights that induce a weighting on the set of spanning trees.  We 
show that the weighted, directed spanning trees 
(often called arborescences) of any planar
graph $G$ can be put into a one-to-one weight-preserving 
correspondence with the perfect matchings of a related
planar graph $\H$.

One special case of this result is a bijection 
between perfect matchings of the hexagonal 
honeycomb lattice and directed spanning trees
of a triangular lattice.  Another special case gives a correspondence 
between perfect matchings of the ``square-octagon'' 
lattice and directed weighted spanning trees on a directed weighted version of the cartesian lattice.  

In conjunction 
with results of \cite{kenyon:dimers}, our main theorem allows us to
compute the measures of all cylinder events for random spanning trees
on any (directed, weighted) planar graph.  Conversely, in cases where
the perfect matching model arises from a tree model, Wilson's
algorithm allows us to quickly generate random samples of perfect
matchings.
\end{abstract}

\section{Introduction}

\cite{temperley:tree-q} observed that asymptotically the
$m\times n$ rectangular grid has about as many spanning trees as the
$2m\times 2n$ rectangular grid has perfect matchings (dimer
coverings).  Soon afterwards he found a bijection between 
spanning trees of the
$m\times n$ grid and perfect matchings in the $(2m+1)\times (2n+1)$
rectangular grid with a corner removed \citep{temperley:tree}.  
The second author of the present article and,
independently, \cite{burton-pemantle:tree}
generalized this bijection to map spanning trees of general
(undirected unweighted) plane graphs to perfect matchings of a related
graph.  Here we extend this bijection to the directed weighted case.  

This generalized bijection can be viewed as a way of
``reducing'' planar spanning tree systems to planar dimer systems 
(though not vice versa in general): 
for any graph whose spanning trees we are interested
in, there is a related graph whose dimer coverings are in a natural
one-to-one weight-preserving correspondence with the 
spanning trees of the original graph.  
Thus properties of spanning trees on any planar graph
can be studied by considering the related dimer system.  
However, only certain dimer systems 
are related to spanning tree systems in the aforementioned way.
Two important examples are perfect matchings of finite subgraphs
of the hexagonal honeycomb lattice (combinatorially equivalent
to ``lozenge'' tilings of finite regions; see e.g.\ \cite{kuperberg:lozenge}) 
and perfect matchings of finite subgraphs of the ``square-octagon'' 
lattice.  
Both of these dimer models are in bijection with weighted, 
directed spanning trees on associated graphs.

There are a number of important applications of our bijection.
Some questions about spanning tree models do not seem to be
amenable to direct analysis, but can be approached if one
first translates the problem into one involving the associated
dimer model and then makes use of tools available in that context.
Conversely, some problems involving dimers are most easily handled 
if one converts them into problems involving spanning trees (though
this can be done only for a limited class of dimer models).
We now describe these applications in greater detail.

One example of a spanning tree property that is easy to study after
reducing the problem to that of dimers is the computation of certain
probabilities, such as
the probability that a directed edge $e_1$ is in the tree
and the directed dual edge $e_2$ is in the dual tree.
(For a definition of dual tree, see \sref{gen-temperley}.)
The presence or absence of the dual edge $e_2$ in the dual tree
is not a local event with respect to the (primal) tree model;
that is, the event is not determined by the presence or absence
of a fixed set of edges in the primal tree.
(The fact that $e_2$ is an {\it oriented\/} edge is crucial here.)
On the other hand, the event in question
is a local event in the associated matching process,
since the matching directly incorporates both primal and dual directed trees.
The probabilities of local events in either the tree or matching model are
easy to compute (\cite{burton-pemantle:tree}, \cite{kenyon:dimers}),
but events of the above type are harder if not
impossible to compute from the point of view of the tree only
\citep{burton-pemantle:tree}.

Another spanning tree property that can be studied via dimers is the
number of times that the path connecting two points in a spanning tree
winds around the two points.  In \sref{winding} we relate these
winding numbers to height functions in the dimer model; the first
author has shown in \cite{kenyon:dimers} how to compute properties of
these height functions (and consequently the corresponding winding
numbers) such as the variance.

Dimer systems can also be studied via trees, if the given dimer system
has a spanning tree model associated with it.  For instance, one can
sometimes enumerate the dimer coverings of a graph by counting the
number of spanning trees in the associated tree model.  In \sref{formula} 
we show a variety of such graphs, together with exact formulas for the number
of dimer coverings, where the easiest (or only) way we know to obtain
these formulas is to count spanning trees. In the dimer model on a bounded
region, the boundary can have an important (long-range) 
effect on the number of configurations \citep{cohn-kenyon-propp}.
In this case the regions which arise from the associated spanning
tree process give the most ``natural'' boundary conditions for the dimer model,
in the sense that the boundary has the least long-range influence
\citep{kenyon:confinv}.

Another case where a spanning tree model is useful for studying the
associated dimer model is in the generation of random samples.
Wilson's algorithm \citep{propp-wilson:tree} can be used to generate
random spanning trees quickly --- the expected running time is given
by the sum of two mean hitting times.  For the lattice of octagons and
squares, the expected running time is {\it linear\/} in the number of vertices.
For the usual lattice of squares, when a rectangular region has
moderate aspect ratio, the running time is nearly linear, but with a
logarithmic correction factor.

Finally, \cite{burton-pemantle:tree} prove that the uniform measure
on spanning trees of the $n\times n$ square grid converges as $n\to\infty$
to the unique translation-invariant measure of maximal entropy on 
the set of spanning forests with no finite component.
Consequently the associated dimer model (on $\Z^2$)
has a unique translation-invariant measure of maximal entropy.
We do not know how to prove this directly from the dimer model itself,
or in any other dimer model except those arising from our construction
via a bijection with {\em undirected} (but possibly weighted) spanning trees.

We remark that other combinatorial systems that can also be reduced to
dimer systems in a similarly simple way include the Ising model on
planar graphs \citep{fisher:ising-dimer}
and systems of non-intersecting lattice paths \citep{lindstrom,gessel-viennot}.

In \sref{gen-temperley} we prove the generalized version of
Temperley's bijection.  
In the two succeeding sections 
(\S\S~\ref{sec:hexagonal}-\ref{sec:square-octagon})
we illustrate the bijection with two examples:
In \sref{hexagonal} we exhibit a bijection between directed spanning trees
on the triangular lattice and
perfect matchings of the hexagonal honeycomb lattice.  Our 
bijection cannot be applied directly to matchings on the
square-octagon lattice, but in \sref{square-octagon} we show how to
locally transform this lattice so that the bijection can be applied.
This transformation enables the rapid generation of random dimer configurations of the square-octagon lattice.
Then in
\sref{winding} we show how the winding number of arcs in a 
spanning tree can be related to the height function on the corresponding 
perfect matching.
In \sref{formula} we use our generalized bijection to compute the exact
number of perfect matchings of some ``locally symmetric'' finite planar 
graphs, that is, graphs that arise as finite induced subgraphs of highly 
symmetric infinite planar graphs.  Lastly, in \sref{open} we give
some open problems.

\section{Generalized Temperley Bijection} \label{sec:gen-temperley}

Let $G$ be a finite connected directed graph embedded in the plane,
with multiple edges and self-loops allowed.  In general the edges of $G$ 
will be weighted, that is, each directed edge from vertex $u$ to vertex $v$
has a nonnegative weight assigned to it, which need not be the same as the
weight of other directed edges from $u$ to $v$ or from $v$ to $u$.
Undirected graphs can be fit into our framework by thinking of 
each undirected edge as two directed edges, one in each direction,
embedded in the plane so as to coincide.
(We will discuss issues related to choice of embedding
at the end of this section.)
Unweighted graphs can be fit into our framework by assigning
each edge weight 1.

By a {\bf directed spanning tree} (or {\bf arborescence}) $T$ of $G$ we mean 
a connected, contractible union of (directed) edges such that 
each vertex of $G$ except one has exactly one outgoing edge in $T$.  
Note that the exceptional vertex necessarily has no outgoing edges in $T$;
it is called the {\bf root} of $T$.  
We define the {\bf weight} of such a tree $T$ to be
the product of the weights of its edges.

\begin{figure}[phtb]
\centerline{\epsfig{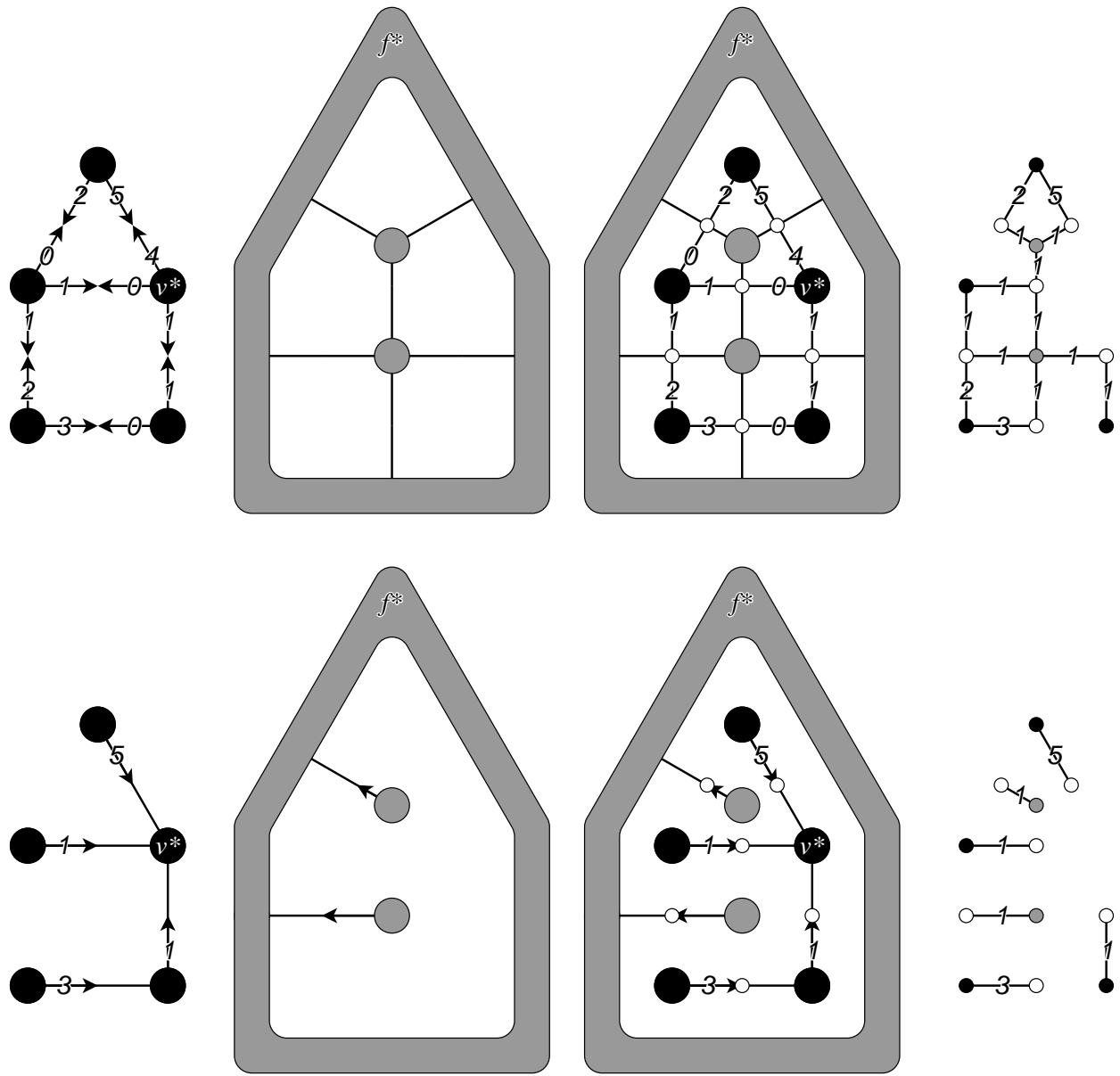}}
\caption{
Illustration of the generalized Temperley bijection.}
\label{fig:fig1}
\end{figure}

We will make a new weighted graph $\cH(G)$ based on $G$, 
as shown in the top half of \fref{fig1}.
$G$ is shown in the top left-most panel.
$G^\perp$, the dual graph of $G$ (second panel),
has vertices, edges, and faces of $G^\perp$
corresponding to faces, edges, and vertices of $G$, respectively
(including a vertex, here marked $f^*$, 
that corresponds to the unbounded, external face of $G$,
and is represented in ``extended form'',
i.e., as a spread-out region rather than a small dot).
We can embed $G$ and $G^\perp$ simultaneously in the plane,
such that an edge $e$ of $G$ crosses
the corresponding dual edge $e^\perp$ of $G^\perp$ exactly once
and crosses no other edge of $G^\perp$.
If we introduce a new vertex at each such crossing,
we get the graph shown in the third panel.
This is the graph $\cH(G)$.
Pictorially, we may derive $\cH(G)$ from $G$
by adding a new node on each edge $e$
and a new node on each face $f$
and joining them by a new edge
if $e$ is part of the boundary of $f$.
To avoid confusion, we will say that $\cH(G)$ has nodes
and links whereas $G$ has vertices and edges.

Here is an alternative, direct definition of $\cH(G)$
that does not go by way of the dual graph.
Put $V=$ the set of vertices of $G$,
$E=$ the set of edges, $F=$ the set of faces
(including the unbounded face).
Define $\cH(G)$ as the weighted undirected graph
with a node $\overline{v}$ corresponding to each vertex $v$ of $G$,
a node $\overline{e}$ corresponding to each edge $e$ of $G$,
and a node $\overline{f}$ corresponding to each face $f$ of $G$,
with a link joining two nodes in $\cH(G)$
if the corresponding structures in $G$ are either
an edge and one of its endpoints
{\it or\/} an edge and one of the faces it bounds.

The weight of a link between a vertex-node $\overline{v}$ and an edge-node 
$\overline{e}$ (where $v$ is an endpoint of $e$ in $G$) is the weight of 
edge $e$ in $G$ directed away from $v$.  
The weight of a link between a face-node $\overline{f}$ and an edge-node 
$\overline{e}$ (where $e$ bounds face $f$ in $G$) is always $1$.

A {\bf perfect matching} of a graph $H$ is a collection of edges $M$
such that each vertex is a vertex of exactly one edge of $M$.
The {\bf weight} of a perfect matching
is the product of the weights of its edges
(1 by default in the unweighted case).

In the case of both trees and matchings,
the weighting gives rise to a probability distribution
on the objects in question,
in which the probability of any particular object
(tree or matching)
is proportional to its weight.

Let $v^*$ be a vertex of $G$
and $f^*$ a face of $G$,
and let $\H = \H(v^*, f^*)$ be the induced subgraph of $\cH(G)$ 
obtained by deleting the nodes $\overline{v^*}, \overline{f^*}$ 
(along with all incident edges in $\cH(G)$),
as shown in the fourth panel of the top half of \fref{fig1}).
Since by Euler's formula
$(|V|-1)+(|F|-1)=|E|$,
$\H(v^*,f^*)$ is a balanced bipartite graph, so it may have perfect matchings.
(For a nice tree-based proof of Euler's formula, see \cite[page 57]{aigner-ziegler:book}.)

\begin{theorem}\label{thm:main}
If $v^*$ is incident with $f^*$, then there is a
weight-preserving bijection between spanning trees of $G$ rooted at
$v^*$ and perfect matchings of $\H(v^*, f^*)$.  If $v^*$ is not
incident with $f^*$, there remains a weight-preserving injection from
the spanning trees of $G$ rooted at $v^*$ 
to the perfect matchings of $\H(v^*, f^*)$.
\end{theorem}

This theorem, along with its proof, is a generalization of a result of
\cite{temperley:tree} which is discussed in problem 4.30 of
\cite[pages 34, 104, 243--244]{lovasz:problems}.  The unweighted
undirected generalization was independently discovered by
\cite{burton-pemantle:tree}, who applied it to infinite graphs,
and also by F.~Y.~Wu, who included it in lecture notes for a course.

Note that in the special case when we take all weights of $G$ to be $1$, 
the first part of the theorem
implies that the number of perfect matchings of $\H(v^*, f^*)$
is independent of $v^*$ and $f^*$, provided that $v^*$ and $f^*$
are incident with one another.  

Henceforth, we refer to 
perfect matchings as simply ``matchings,''
and directed spanning trees of $G$ rooted at $v^*$ as simply 
``spanning trees''
or occasionally just ``trees.''
\medskip

{\sc Proof of theorem:}
It will be enough to exhibit a weight-preserving injective mapping
from the set of spanning trees of $G$
into the set of matchings of $\H(v^*,f^*)$,
and to show that when $v^*$ is incident with $f^*$,
every matching of $\H(v^*,f^*)$ arises from a spanning tree of $G$.

Given a spanning tree $T$ of $G$ rooted at $V^*$,
the set of edges of $G^\perp$ that do not cross edges of $T$
form a spanning tree of $G^\perp$, called the {\bf dual tree} 
and here denoted by $T^\perp$.
Orient the edges of $T^\perp$ so that
they point towards $f^*$.
Then a matching $M$ of $\H(v^*,f^*)$ can be obtained
as shown in the bottom half of \fref{fig1}.
Specifically, for each $v \in V$,
pair $\overline{v}$ with the unique $\overline{e}$
such that $v$ is an endpoint of $e$
and $e$ is pointing away from $v$ in the orientation of $T$,
and for each $f \in F$,
pair $\overline{f}$ with the unique $\overline{e}$
such that $e$ bounds $f$
and $e^\perp$ is pointing away from $f$ in the orientation of $T^\perp$.
The left panel shows the tree $T$;
the second panel shows the dual tree $T^\perp$;
the third panel shows both trees;
and the fourth panel shows the matching $M$,
which has the same weight as $T$.

To verify that this construction always gives a matching $M$ of $\H(v^*,f^*)$,
it suffices to show that no edge-node $\overline{e}$ is paired twice.
But this could only happen if we had $e \in T$ and $e^\perp \in T^\perp$,
contradicting the definition of a dual tree.

{}From the matching $M$ we can easily recover $T$ as the set of edges
$e$ such that $\overline{e}$ is paired with a vertex-node in $\H(v^*,f^*)$ under
the matching $M$.  Hence the mapping $T \mapsto M$ is injective.

Now suppose $v^*$ is incident with $f^*$,
and let $M$ be a matching of $\H(v^*,f^*)$.
Let $\widetilde{T}$ be the set of edges $e$ of $G$
such that $\overline{e}$ is paired with a vertex-node by $M$.
To complete the proof of the theorem,
we must show that $\widetilde{T}$ is a spanning tree.
Note that $\widetilde{T}$ has $|V|-1$ edges,
so it suffices to prove that $\widetilde{T}$ is acyclic.

Suppose $\widetilde{T}$ contained a cycle $C$, say of length $n$.
$C$ divides the plane into two (open) regions,
one of which contains both $v^*$ and $f^*$
and the other of which contains neither.
We claim that each part contains
an odd number of nodes of $\cH(G)$
and hence an odd number of nodes of the subgraph $\H(v^*,f^*)$ as well.
For, suppose we modify $G$ by replacing
either of the two regions by a single face.
By Euler's formula,
the number of vertices, edges, and faces
in the resulting graph must be even.
Since there are an even number of these elements
on the cycle $C$
($n$ vertices and $n$ edges)
and an odd number in the modified region (1 face),
the unmodified region must have
an odd number of elements as well.

Since the edges of $C$ disconnect $\H(v^*,f^*)$ into parts
lying in the two regions,
$M$ must match each region within itself.
But this is impossible,
since each region has been shown
to contain an odd number of nodes of $\H(v^*,f^*)$.
This completes the proof of the theorem.
$\Cox$

As was remarked earlier,
the theorem implies that
when $v'$ is incident with $f'$
and $v''$ is incident with $f''$,
the matchings $M'$ of $\H(v',f')$ are equinumerous with
the matchings $M''$ of $\H(v'',f'')$;
in fact, the proof of the theorem
provides a bijection between the two sets of matchings.
This bijection can be understood without reference to spanning trees,
as a process of ``sliding edges.''
Specifically, one iteratively defines a chain
$v'' = v_0,\, e_0,\, v_1,\, e_1,\, v_2, ...$
such that, for all $i$,
$\overline{e_i}$ is the node that $M'$ pairs with $\overline{v_i}$
and $v_{i+1}$ is the vertex of $e_i$ that is distinct from $v_i$.
This chain cannot repeat any vertices,
since any closed loop would encircle an odd number of nodes
(see the preceding proof),
so it must terminate by arriving at $v'$ after some number of steps.
That is, the chain must be of the form
$$v'' = v_0,\, e_0,\, v_1,\, e_1,\, v_2, ..., e_{r-1},\, v_r = v'$$
for some $r$.
Once one has found such a chain,
one modifies the matching $M'$
by pairing $\overline{e_i}$ with $\overline{v_{i+1}}$ 
instead of $\overline{v_i}$.
One then does the same with a chain of dual-edges
joining the faces $f''$ and $f'$,
obtaining the desired matching $M''$.

We also remark that in addition to one's having a choice of which
vertex-node and face-node to delete, one often has a choice of how
to embed a graph in the plane in the first place.  For instance,
in the case where $G$ has a single edge from $u$ to $v$ and a single
edge from $v$ to $u$, we allowed the two edges to be embedded so as 
to coincide.  What if we had required the embedding to be proper, so
that the two edges could meet only at their endpoints?  Then one
would get a slightly enlarged graph $\cH(G)$ in which a single edge-node
in the original $\cH(G)$ was replaced by two edge-nodes and a face-node
in between (corresponding to the digon bounded by the two edges).
It is easy to see in this case that perfect matchings of the first
$\H(G;v^*,f^*)$ are in bijection with perfect matchings of the second
$\H(G;v^*,f^*)$.
When there are multiple directed edges in each direction, the number of
possible embeddings increases rapidly; but our main bijection theorem
guarantees that the number of matchings of $\H(G;v^*,f^*)$ is insensitive
to the choice of embedding.

Moreover, having several directed edges from $v$ to $w$ is in a certain
sense equivalent to having a single edge from $v$ to $w$ whose weight is
the sum of the weights of those directed edges.  It is not true that
the spanning trees of the former graph are in bijection with those of
the latter graph; however, there is an obvious mapping from the former
to the latter, and this correspondence is weight-preserving, in the
sense that the weight of a spanning tree of the smaller graph is the
sum of the weights of the spanning trees in the larger graph to which
it corresponds.  It follows that the sum of the weights of all the
spanning trees is the same for both graphs.

Given a graph $H$, it can be an amusing problem to find a directed
graph $G$ such that $\H(G;v^*,f^*)=H$.  We leave it to the reader to show that
this cannot be done with the square-octagon lattice of \fref{square-octagon}.
(That is, there is a finite subgraph of the lattice, such that any
subregion $H$ of the square-octagon lattice containing this subgraph will
fail to be of the form $\H(G;v^*,f^*)$.)

\section{The Hexagonal Lattice} \label{sec:hexagonal}

In this section we illustrate the technique of \tref{main} by
giving a bijection between spanning trees of a directed graph and
matchings in the hexagonal (honeycomb) lattice.  

\begin{figure}[phtb]
\centerline{\epsfig{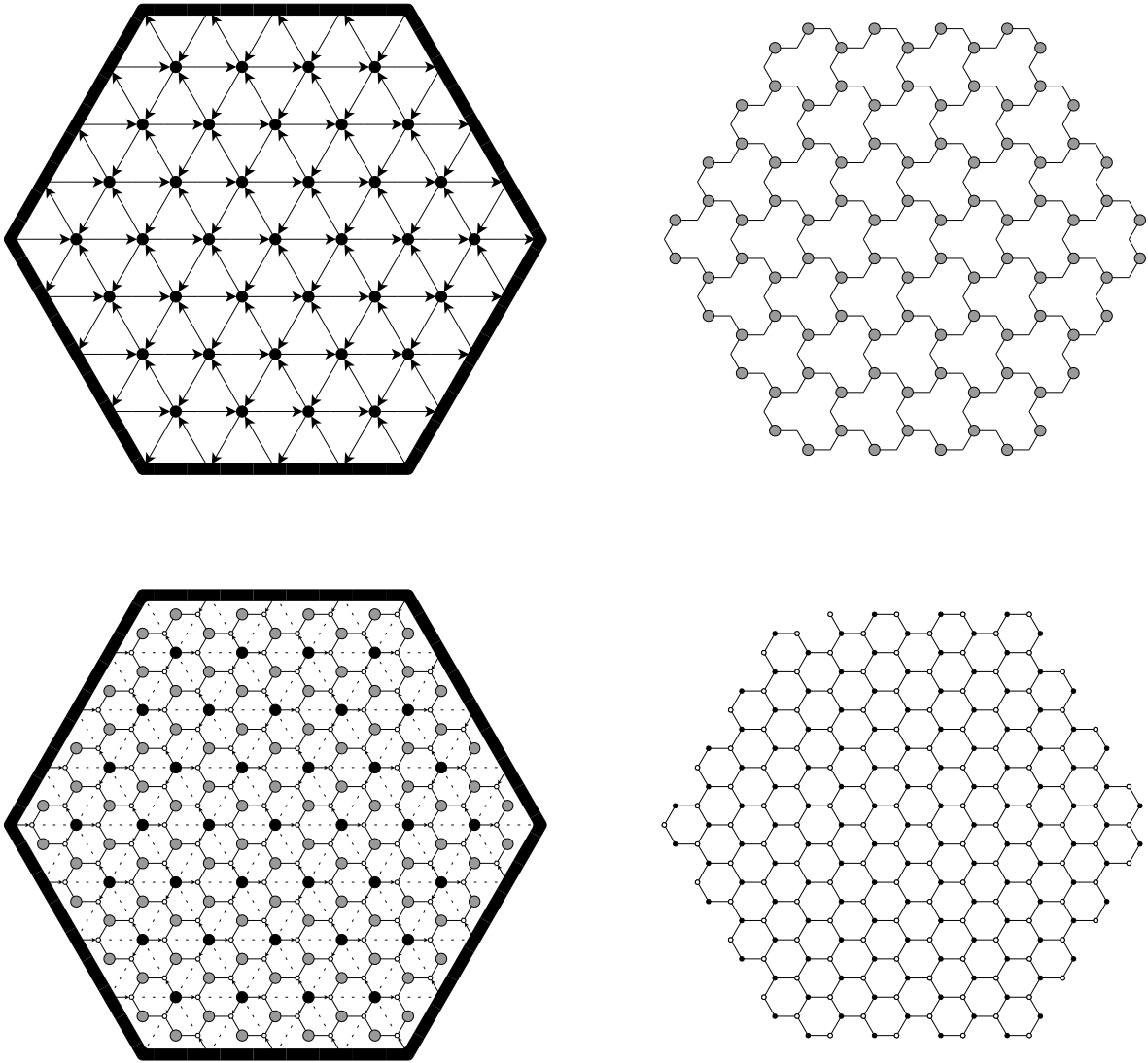}}
\caption{Generalized Temperley bijection for the hexagonal lattice.}
\label{fig:hextree}
\end{figure}

Panel (a) of \fref{hextree}
contains the plane graph $G$, a directed triangular lattice.
(Here and throughout the rest of the article, the upper-left,
upper-right, lower-left, and lower-right panels of a four-panel
figure will be denoted by (a), (b), (c), and (d), respectively.)
$G$ contains an ``outer vertex'' 
which is represented in extended form, in this case drawn as
a large hexagon.  In the examples throughout the rest of the article,
either $G$ or $G^\perp$ (or both) will have an outer vertex that is
drawn in extended form.  Panel (b) shows the dual of $G$, a hexagonal
lattice.  The edges in
panel (b) have been drawn bent slightly so that the union of $G$ and
its dual (panel (c)) can be recognized as a subset of the hexagonal
lattice.  The dotted edges in panel (c) have weight zero, and may be
omitted; they are shown only to highlight the connection with panel (a).  
The graph $\cH(G)$ can be read off from panel (c);
it is a hexagonal lattice with about three times as many hexagons
as $G^\perp$.
Panel (d) shows the graph $\H(v^*,f^*)$,
which is obtained from $\cH(G)$ by removing $v^*$ (the outer vertex of $G$) 
and $f^*$ (the leftmost vertex at the top of $G^\perp$).
We shall apply this correspondence in \S~\ref{ssec:expl-hex}.

\section{The Square-Octagon Lattice} \label{sec:square-octagon}

Here we illustrate a less direct application of \tref{main},
and give a bijection between perfect matchings of certain
planar graphs and spanning trees on an associated graph. 
Consider perfect matchings on the square-octagon lattice, an excerpt 
of which is shown in \fref{square-octagon}.
This graph does not arise as $\cH(G)$ for any graph $G$, so 
\tref{main} does not apply immediately.
\begin{figure}[htb]
\centerline{\epsfig{figure=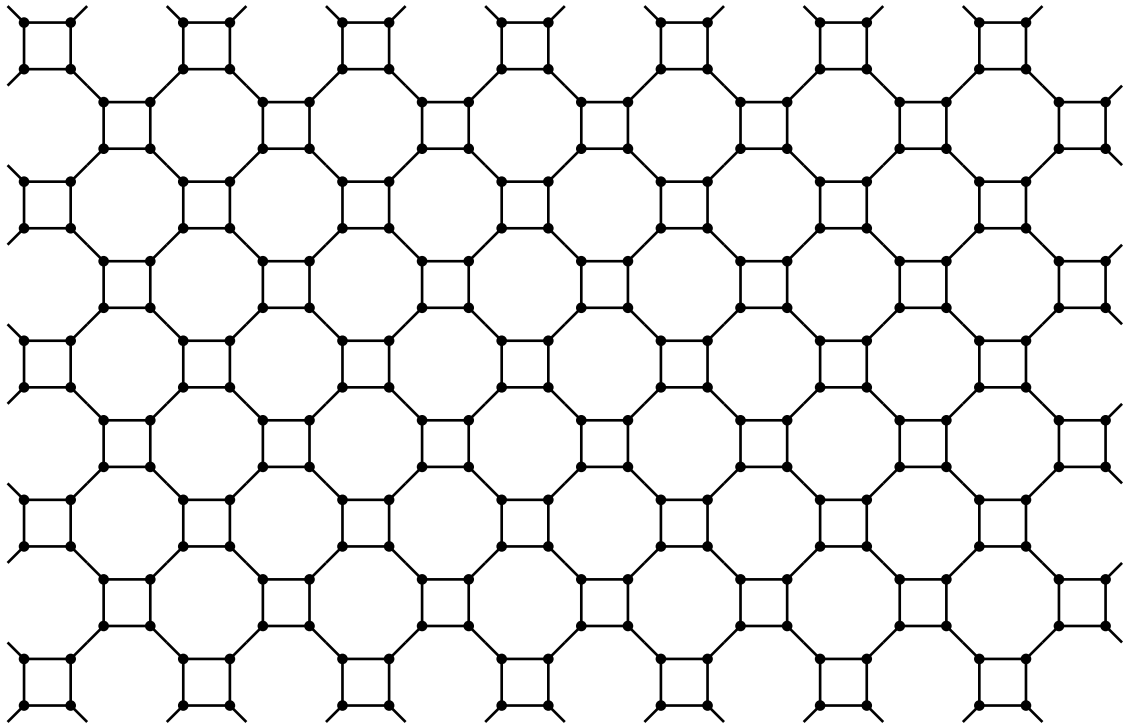}}
\caption{A portion of the square-octagon lattice.}
\label{fig:square-octagon}
\end{figure}
Nonetheless, it is possible to
generalize the bijection to 
apply to
this lattice.  To do it we need to apply two
transformations to the lattice.  The first transformation is called
``urban renewal,'' a term coined by the second author,
who learned of the method from
Greg Kuperberg.  In the second transformation, we adjust the edge
weights.  At that point, if a square-octagon region has suitable boundary 
conditions,
the transformed graph can be expressed as $\cH(G)$ for some graph $G$.

\subsection{Urban renewal}

Tricks such as urban renewal have
been used by researchers in the statistical mechanics literature for
decades, but since understanding it is essential for what follows, a
description is included here
of the special case of urban renewal that we will need. 
One views the
square-octagon lattice as a set of cities (the squares) that
communicate with one another via the edges that separate octagons.
Now the graph of cities 
(with each city being thought of as adjacent to the four closest cities)
is itself bipartite, so we may say that every
city is either rich or poor, with every poor city having four rich
neighbors and vice versa.  The process of urban renewal on a poor city
merges each of its four vertices with its four neighboring vertices,
and then changes the weights of the edges of the poor city 
from 1 to $1/2$, as shown in \fref{urban}.
We will show that 
the sum of the weights of the matchings in the ``before'' graph
is twice the sum of the weights of the matchings in the ``after'' graph.
We will do this by associating with each matching in the before graph
one or two matchings in the after graph,
and vice versa.
More precisely,
we divide the set of matchings in the before graph into equivalence classes
of size 1 or 2,
and likewise with the set of matchings of the after graph,
and we create a bijection between these equivalence classes
so that the weight of each class in the before graph
(that is, the sum of the weights of the matchings
that constitute that class)
is twice the weight of the associated class in the after graph.

\begin{figure}[htb]
\centerline{\epsfig{figure=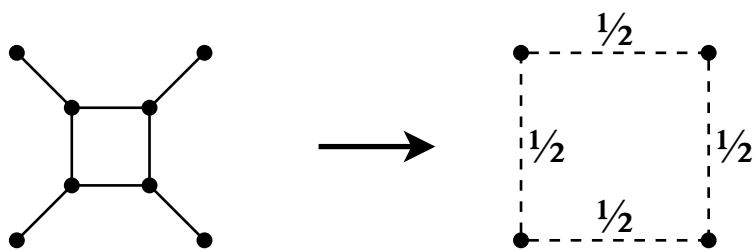}}
\caption{Urban renewal.  The poor city is the inner square in the left graph, 
and is connected to the rest of the graph (not shown) via only four vertices 
at the corners, some of which may actually be absent.  The city and its 
connections are replaced with weight $1/2$ edges, shown as dashed lines.
All other edges have weight 1.}
\label{fig:urban}
\end{figure}

Matchings in the ``before'' graph get mapped via urban renewal to matchings 
in the ``after'' graph by deleting the four vertices of the poor city
and its incident edges, and then pairing up any resulting unpaired
vertices.  Prior to urban renewal, every matching will match $k$ of
the poor city's vertices with the rest of the graph, with $k$ equal to
0, 2, or 4; if $k=2$, then these vertices are adjacent.  If $k=0$,
then since the city has two possible matchings, a pair of matchings in
the ``before'' graph get mapped to one matching (of half their
combined weight) in the ``after'' graph.  If $k=2$ (two of the poor
city's vertices match to each other and two match outward), then the
matching in the before graph gets mapped to a matching in the after
graph that uses one weight-$1/2$ edge.  The matchings with $k=4$ get
mapped to a pair of matchings in the after graph, each using two
weight-$1/2$ edges.  Thus urban renewal on a poor city will reduce the
weighted sum of matchings by a factor of $1/2$.  (If one is trying to
generate random matchings rather than merely count them, then, given a
random bit, a random matching in the before graph is readily
transformed into a random matching in the after graph, and conversely,
given a random bit, a random matching in the after graph is readily
transformed into a random matching in the before graph.)

The preceding discussion applies to cities in the interior of
a finite subgraph of the infinite square-octagon grid.
Along the boundaries, some of the poor cities may not have four
neighbors, but urban renewal can still be done.  
One way to see this
is to adjoin a pair of connected vertices to the graph for each
missing poor city's neighbor, and connect one of these vertices to the
poor city.  This operation won't affect the number of matchings or
their weights, and after urban renewal, the pair may be
deleted again, again without affecting the matchings --- so if some of
the poor city's vertices don't have neighbors, these vertices are
deleted by urban renewal.  

Doing urban renewal on each of the poor
cities in the square-octagon lattice will yield the more familiar
Cartesian lattice.

\subsection{Weighted directed spanning trees}

Consider the finite square-octagon graph shown in \fref{diab-lm}.  
It has 3 octagonal bumps on the left, and
four on top, so by convention let's call it a region of order $3,4$.
(In a region of order $L,M$, there are $2LM$ octagons.)  An octagonal
column and octagonal row meet at a unique square; these will be the
rich cities.  The rich cities have been labeled by their coordinates
to enhance clarity.  The other $(L+1)(M+1)$ squares will be the poor
cities, and we will do urban renewal on them as shown in
\fref{diab-lm}.  We will compute the weighted number of
matchings of the resulting graph, and multiply by $2^{(L+1)(M+1)}$.

\begin{figure}[tbp]
\centerline{\epsfig{figure=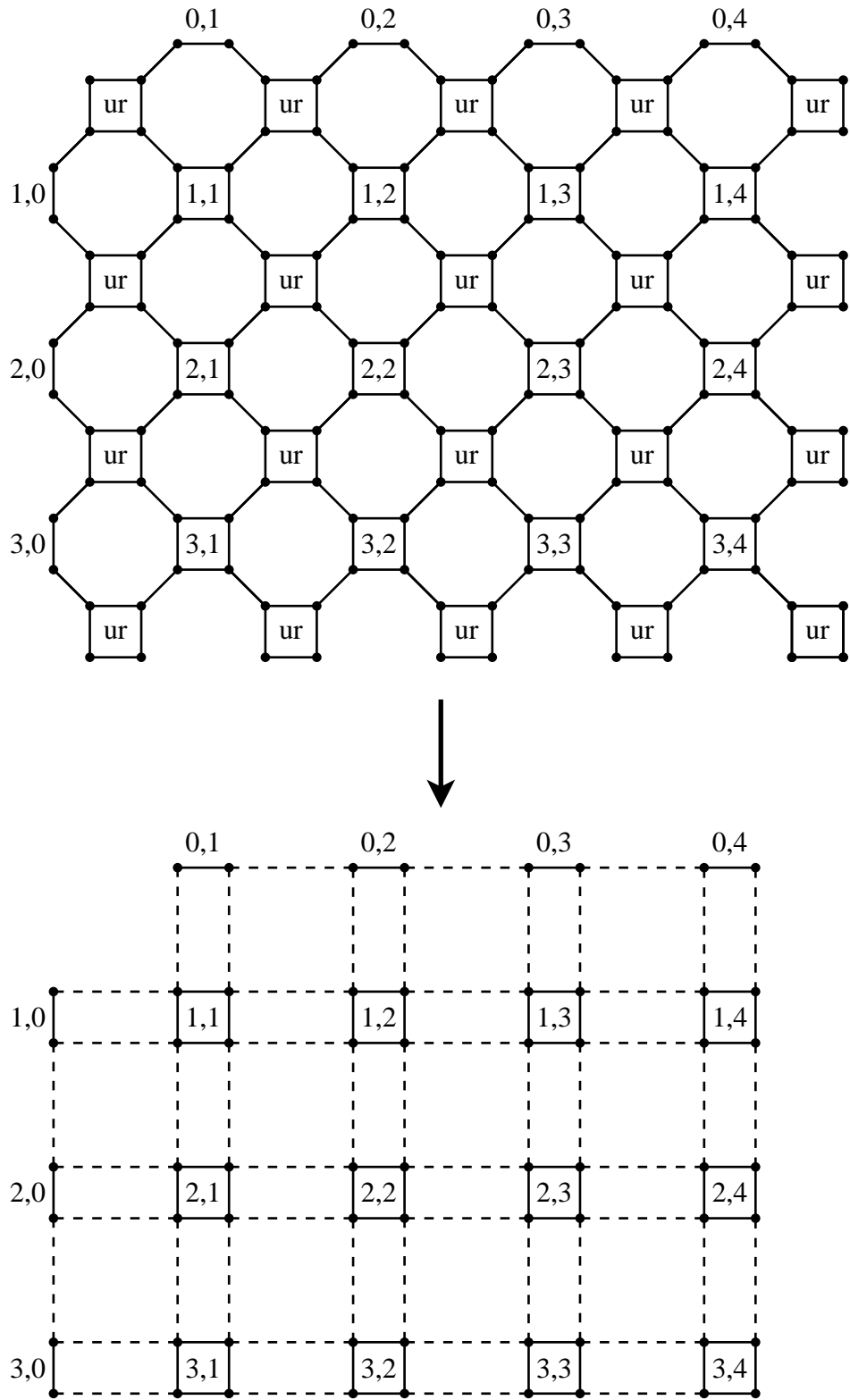}}
\caption{Region of order $3,4$ before and after urban renewal.  The poor 
cities on which urban renewal is done are labeled with ``ur,'' the rich 
cities are labeled with their coordinates.  Dashed edges have weight $1/2$.}
\label{fig:diab-lm}
\end{figure}

Now for any vertex, we may re-weight all the edges incident to the
vertex (multiplying them all by the same constant)
without affecting the probability distribution on matchings:
this has the effect of multiplying the weight of any matching by
that same constant.
Re-weight the edges as follows: 
for the rich city (complete or incomplete) with coordinates $i,j$
($0 \leq i \leq M$, $0 \leq j \leq N$),
multiply the weights of the edges incident to the top left and
lower right corners by $2^{-i-j}$, the other two corners by $2^{i+j}$.

Edges that are internal to the rich cities remain weighted by 
$2^{i+j} 2^{-i-j} = 1$.  
The
long edges come in pairs.  The lower or right edge of the pair gets
its weight doubled, to become 1, while the upper or left edges of the
pair gets its weight halved to become 1/4.

The next thing we need to do is interpret this graph as a plane graph
and its dual (see \fref{diabtree}).  The upper left vertices
of the small squares represent vertices, the lower right vertices
represent faces, and the other two vertices represent edges.  The
result is the graph 
$F$
shown in \fref{diabtree}, which has
$LM+1$ vertices --- $LM$ of them on a grid, and one ``outer vertex''
(not in the original graph) that all the open edges connect to.  A
random spanning tree on the vertices of this graph rooted at the outer
vertex determines a dual tree on the faces of this graph, rooted at
the upper left face, 
and the two together determine
a matching of the
graph in \fref{diab-lm}.  The weight of the matching equals the weight of the
primal tree, since the re-weighting left every dual tree with weight
one.

\begin{figure}[tbp]
\centerline{\epsfig{figure=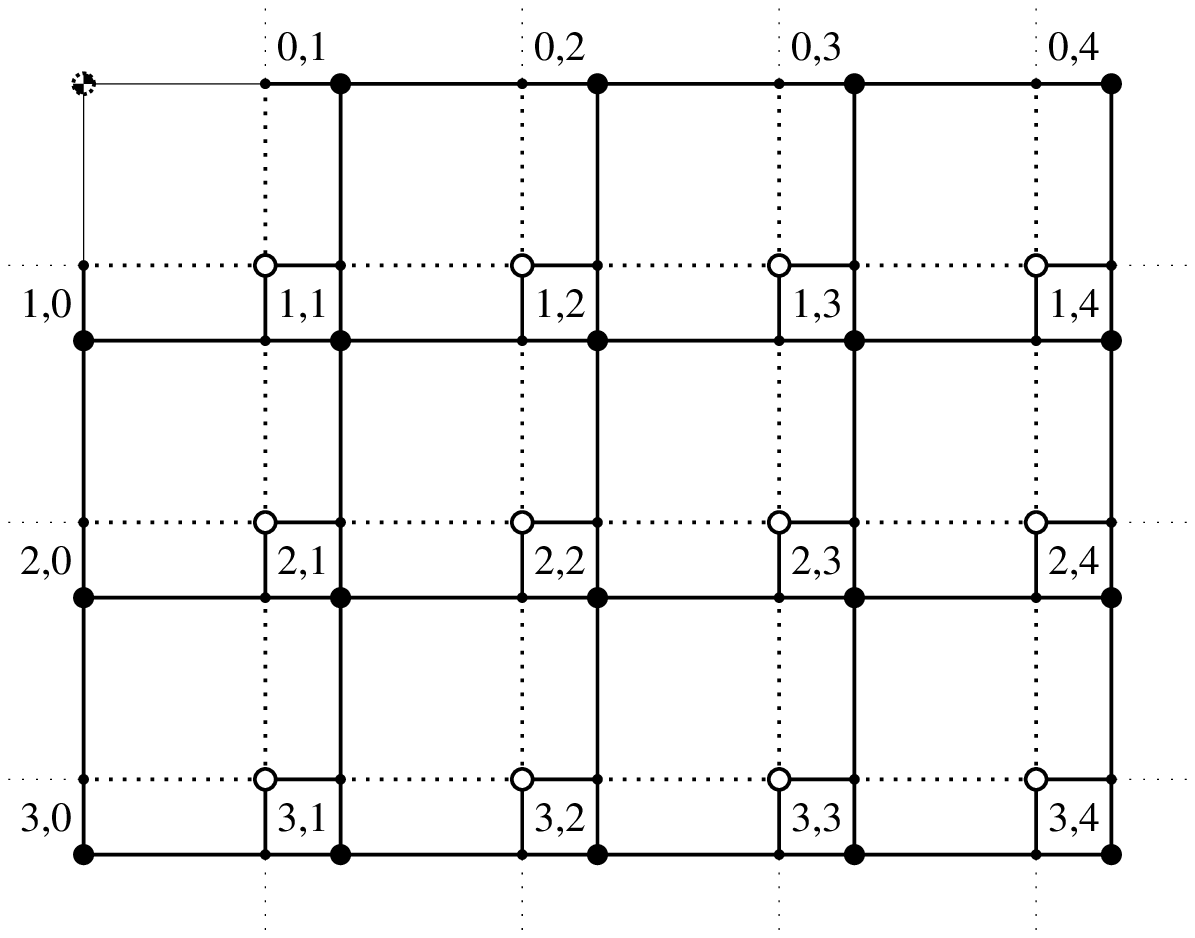}}
\caption{The region from Figure~\protect\ref{fig:diab-lm} after re-weighting. 
Dashed edges have weight $1/4$.  
A distinguished vertex and distinguished 
outer face have been adjoined to give a graph $H$ which nicely decomposes 
into a graph $G$ and a weighted directed graph $F$.}
\label{fig:diabtree}
\end{figure}

\subsection{Random generation in linear expected time} 
\label{ssec:linear}

Using the loop-erased random-walk spanning tree generator
\citep{propp-wilson:tree}, and the bijection derived above
between spanning trees of a weighted graph and matchings of the
square-octagon regions, 
we can sample random matchings in linear time.
The tree generator builds the tree by doing a sequence of
loop-erased random walks on the underlying graph.
(From any vertex $v$, the probability of moving to a particular
neighboring vertex $w$ is proportional to the weight of the edge
from $v$ to $w$;
this determines a random walk on the graph.
For details on loop-erasure, see 
\citep{propp-wilson:tree}.)
It has been shown that
the expected running time (or rather, number of random-walk steps) of
the tree algorithm is given precisely by $$\sum_v E\left[\parbox{2.5in}{\#
times a that random walk started at $v$ visits $v$ before hitting the
root}\right].$$ For our random walk, the moves are right or down with
probability 4/10 and up or left with probability 1/10, since in the
face graph the links going to the left or up have 1/4 the weight of
the other links.  For large graphs, the random walk drifts to the
right and down, so we consider this biased random walk on $\Z^2$.
Starting at the origin, with probability 1 the origin is visited
finitely many times.  Let $R$ be the expected number of times the
random walk returns to its starting location, counting the ``return''
at time 0, before drifting off to infinity.  The first expression
below for $R$ is not hard to check, and the remaining two can
be found in \citep[p.~408]{CRC}:
\begin{align*}
R &= \sum_{k=0}^\infty \binom{2k}{k}^2 (1/5)^{2k} \\
  &= \int_0^1 \int_0^1 \frac{1}{1-(2/5)\cos(\pi x) -(2/5)\cos(\pi y)} 
  \: dx\ dy \\
  &= (2/\pi) K(4/5)
\end{align*}
with $K(k)$ denoting Legendre's complete elliptic integral of the
first kind.
The expected number of steps to create the random spanning tree
is bounded by $R\doteq 1.27025$ times the number of vertices, and the
remaining steps are readily done deterministically in linear time.

\section{Height Functions and Winding Numbers}\label{sec:winding}
In this section we describe the connection between the
winding number of a spanning tree on a planar graph $G$
and a {\it height function} on its corresponding
matching graph $\cH$.
The result, \tref{wind} below,
answers a question posed to the first author by Itai Benjamini.

\subsection{Height function definition}
We assume that $G$ is connected and that $\cH(G)$ is embedded in the
plane with straight edges.  The straight-line embedding is not
necessary for the definition (see below) but the construction is more
geometric in this case.  Moreover, we assume that $\cH(G)$ is embedded
so that one of $f^*$ or $v^*$ is the outer node (as in \fref{height}),
or else both $f^*$ and $v^*$ are on the outer facet (as in
\fref{2x2height}).

Recall that each facet of $\cH(G)$ is a quadrilateral
containing a node $\bar v$, a node $\bar f$, and two nodes
$\bar e_1,\bar e_2$.  The nodes $\bar v$ and $\bar f$ are opposite
each other. Let $d$ be the diagonal of the quadrilateral facet
directed from $\bar v$ to $\bar f$. 
Let $\arg(d)\in[0,2\pi)$ 
denote the angle of the vector $d$ with respect to the $x$-axis.

Let $M$ be a perfect matching of $\H(v^*,f^*)$ and $T,T^\perp$
respectively the associated spanning tree and its dual.  Let $\D$ be
the set of the diagonals of facets of $\H(G)$.  We will define a
real-valued \textbf{height function} $h\colon\D\to\R$ associated with
matching $M$ (refer to \fref{height} and \fref{2x2height}).

\smallskip

\begin{remark}
In many of our examples in \sref{formula}, one or more vertex or face
nodes are drawn in an extended format.  In these cases, the
``diagonals'' incident to an extended node may be drawn from any point
in the node.  In many situations it is natural to draw more than one
diagonal on a facet if one of the nodes bounding it is drawn in an
extended fashion (as in \fref{2x2height}), and then each of the
diagonals gets its own height.  For instance, to recover the standard
definition of height function for matchings of subgraphs of $\Z^2$, it
is necessary to draw multiple diagonals (see \fref{2x2height}).  It is
thus more natural to view the heights as being defined on the
diagonals of the facets rather than the facets themselves.
\end{remark}

\begin{figure}[phtb]
\centerline{\epsfig{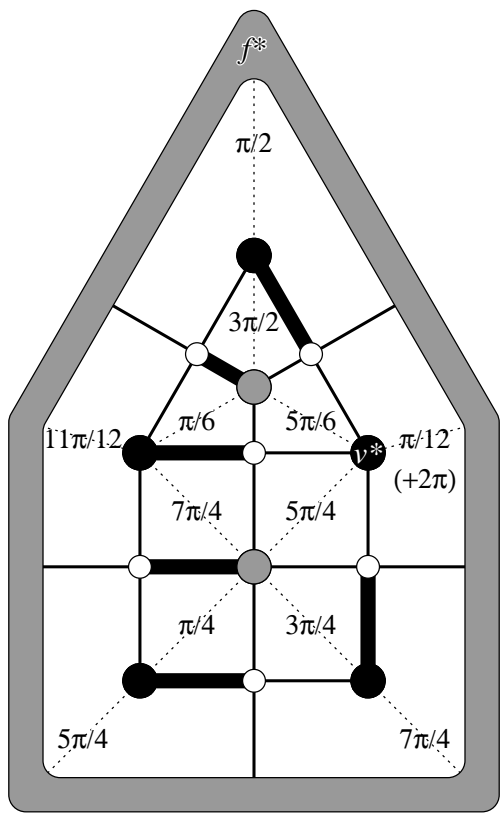}}
\caption{
The height function for the perfect matching in \protect\fref{fig1}.
Since $v^*$ and $f^*$ are unmatched, the height drops by $2\pi$ on
the facet containing $v^*$ and $f^*$.}
\label{fig:height}
\end{figure}

\begin{figure}[phtb]
\centerline{\epsfig{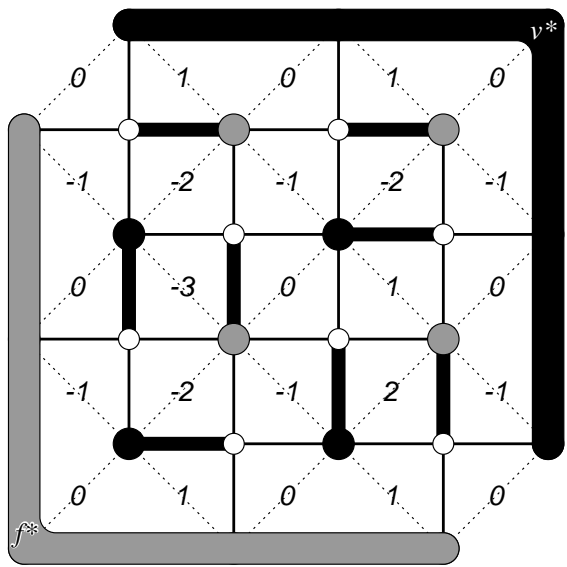}}
\caption{
The height function associated with a perfect matching of another
graph.  For backwards compatibility with previous definitions of the
height function associated with dimers on the square lattice, we have
(1) measured the heights in quarter-turns rather than in radians,
which introduces a scale factor of $\pi/2$, and (2) drawn multiple
diagonals (each with its own height) on facets bounded by an
extended node.}
\label{fig:2x2height}
\end{figure}

We first cut the plane along the links in the perfect matching $M$.
We need every vertex-node and face-node to be at the end of one cut,
and since neither $v^*$ nor $f^*$ is in the matching $M$, we make one
additional cut, from $v^*$ to $f^*$.  If both $v^*$ and $f^*$ border
the outer facet, then we make the cut between them pass through
$\infty$.  If one of $v^*$ or $f^*$ is the outer node, then we view
the outer node as being at $\infty$, so that the cut still goes to
$\infty$.  It is convenient to make this cut split the diagonal from
$v^*$ to $f^*$, so that there are two diagonals in $\D$ from $v^*$ to
$f^*$, one on either side of the cut.

We require the height function $h$ to satisfy the following local
constraint.  Suppose $d$ and $d'$ are two diagonals that share a
vertex node or face node $x$.  Let $\theta$ be the angle required to
rotate $d$ to $d'$ around $x$.  Either the counterclockwise (positive)
rotation or the clockwise (negative) rotation will encounter the cut
containing $x$; we take $\theta$ to be the rotation which avoids the
cut.  The local constraint is $$h(d')=h(d)+\theta.$$

\newpage

\begin{lemma}
Up to a global additive constant, there is a unique height function 
which satisfies the local constraints.  The additive constant can be
chosen so that for each diagonal $d$, $h(d) \equiv\arg(d) \bmod{2\pi}$.
\end{lemma}

{\sc Proof:}
If two diagonals share the same vertex node or face node, then their
height difference is determined by the local constraints.  Suppose
that vertices $v_1$ and $v_2$ of $G$ are connected by an edge, and
that $f$ is a face bounded by this edge.  The height difference
between the diagonal from $v_1$ to $f$ and the diagonal from $v_2$ to
$f$ is determined, and consequently, the local constraints determine
the height difference between any diagonal with $v_1$ as vertex node
and any other diagonal with $v_2$ as vertex node.  Since $G$ is
connected, there can be at most one height function (up to global
additive constant).

We next check that the local constraints do not overconstrain the
height function, i.e.\ that there is a height function satisfying
them.  Between any two diagonals for which there is a local
constraint, one can draw a path connecting the diagonals but which
avoids the cuts.  If the local constraints were inconsistent, then
there would be a closed loop in the plane, which avoids the cuts, such
that the local constraints on the diagonals crossed by the loop are
inconsistent.  Since the cut from $v^*$ to $f^*$ passes through
$\infty$, the interior of this loop does not contain $v^*$ or $f^*$.
Consider such a contradictory loop surrounding a minimal number of
links in the matching $M$.  The loop must cross at least one diagonal
between a vertex node and a face node, and one of these nodes must be
in the interior of the loop.  Call this node $x$.  Since $x$ is not $v^*$
or $f^*$, it is paired with an edge node $y$ in the matching.  Since
the loop avoids cuts, $y$ is also in the interior of the loop.

Since $x$ is a vertex node or a face node, the local constraints on
the diagonals incident to $x$ involve rotations that avoid the cut
from $x$ to $y$, and they are evidently consistent.  Since $y$ is an
edge node, there are four facets of $\H(G)$ incident to $y$, and the
four diagonals in these facets form a quadrilateral containing $y$.
The total height change going {\it clockwise} around $y$ is then the
sum of the interior angles of this quadrilateral, excluding the angle
at the node $x$.  Since the sum of the interior angles of a
quadrilateral is $2\pi$, the total height change around $y$ is $2\pi$
minus the angle at $x$.  Thus the total height change going around the
cut from $x$ to $y$ is constrained to be zero.  Therefore the
contradictory loop can be deformed to exclude $x$ and $y$ from its
interior, contradicting are assumption that it surrounds a minimal
number of links from the matching $M$.  We conclude that there are no
such contradictory loops, and that the height function is well-defined
(up to a global additive constant).

The second statement of the lemma follows by noting that if for some
diagonal $d$, $h(d)\equiv\arg(d) \bmod{2\pi}$, then this relation
holds for the neighboring diagonals as well.
$\Cox$

Note that in the case of a matching of $\Z^2$, this definition of
height function is $2\pi/4$ times the standard definition due to
\cite{thurston:conway} (see \fref{2x2height}).
It is also essentially equivalent to, but more geometrical than,
the definition due to \cite{propp:lattice}.

When $\H(G)$ is embedded in the plane but not with straight edges,
one can still assign to each diagonal $d$ an angle
which is the argument of 
the vector representing the difference in its endpoints.
The orientation of the triple $d,\ell,d'$ is a topological condition
and so does not depend on the fact that the edges of $\H(G)$ are straight.
Thus it is possible to define the height for any embedding of $\H(G)$.

\subsection{Turning and heights}

Let $\gamma$ be a simple path (in topological terms, a ``directed arc'') 
in the spanning tree $T$ from $v'$ to $v''$; 
specifically $\gamma$ is the link path
$$v'=v_0,e_0,v_1,e_1,\ldots,v_{r-1},e_{r-1},v_r=v''.$$
Note that the matching of $\H(v^*,f^*)$ that corresponds to $T$
matches node $v_i$ with node $e_i$ for all $0 \leq i \leq r-1$.
Let $e_r$ be an edge-node (adjacent to $v_r$) 
towards which the path can be continued.

Let $f_1,\ldots f_k$ be the chain of facets of $\H(G)$ 
which share a node with $\gamma$
and lie to the left of $\gamma$, where $f_1$ contains the first link $v_0e_0$
and $f_k$ contains the link $e_{r-1}v_r$. 
Then for $i\in[1,k-1]$
the facets $f_i$ and $f_{i+1}$ are adjacent along
a single link of $\H(G)$, which furthermore is an {\it unmatched} link
(this is a link which shares a node with $\gamma$ but is not in $\gamma$;
therefore it is unmatched). 

The {\em winding number} of $\gamma$ is defined to be the total
angle of the left turns minus the total angle of right turns,
from $v_0$ to $v_r$ (the ``initial'' angle of $\gamma$ is the direction
of $v_0e_0$ and the ``final'' angle of $\gamma$ is the direction
of $v_re_r$).

\begin{theorem}
\label{thm:wind}
The winding number of $\gamma$ is equal to $(h(f_k)-c_k)-(h(f_1)-c_1)$, where
$c_1$ is the counterclockwise angle from the vector $v_0e_0$ to the diagonal 
$d_1$, and $c_k$ is the counterclockwise angle from the vector $v_re_r$ to 
the diagonal $d_k$.
\end{theorem}

Here we may view the heights as being defined on the facets, since
even if facet $f_j$ has more than one diagonal, $h(f_j)-c_j$ has the
same value no matter which diagonal is used.
\medskip

{\sc Proof of theorem:}
Without loss of generality we may assume that the
height at facet $f_1$ is $h(f_1)=\arg(d_1)$. Then $h(f_1)-c_1$
is the angle that the initial direction of $\gamma$ makes with
the $x$-axis.  Let $f_j$ be a facet adjacent to $v_\ell$.
Then similarly $h(f_j)-c_j$ is equal modulo $2\pi$ to the angle that
$v_\ell e_\ell$ makes
with the $x$-axis. 

If $f_j$ is the last facet adjacent to $v_\ell$, so that $f_{j+1}$ is
adjacent
to $v_{\ell+1}$, then the difference $(h(f_{j+1})-c_{j+1})-(h(f_j)-c_j)$
equals the increase in angle from $v_je_j$ to $v_{j+1}e_{j+1}$ (which
is negative at a right turn). The proof follows.
$\Cox$

\section{Explicit Formulas}\label{sec:formula}

Here we use the generalized Temperley bijection to count perfect matchings
of certain finite subgraphs of the infinite square grid and infinite
square-octagon grid, making use of spanning trees.
The enumeration techniques closely follow the
derivation of the exact formula for the number of domino tilings of
the $(2n+1)\times(2n+1)$ square with a corner removed; see 
\cite{lovasz:problems} and \cite{propp:dominoes}.

One motivation for some of these calculations is that they can be used to
compute asymptotic formulas for the number of dimer configurations on
more general regions, via techniques developed in
\cite{kenyon:laplacian}.  For example, using the exact formula for the
triangular or diamond regions in \S\S6.5--6.8, one should be able to
extend the asymptotic formula in
\cite{kenyon:laplacian} to polygonal regions whose boundary
edges have slopes in $\{0,1,-1,\infty\}$.  Using the example in
\ref{ssec:expl-hex} should give rise to a similar asymptotic formula for
regions of the hexagonal lattice.
Other formulas corroborate refinements of the entropy formula that
predict how the the lower-order asymptotics of the number of spanning
trees should reflect the geometry of the boundary of the graph
(see \cite{duplantier-david} and \cite{kenyon:laplacian}).

To enumerate the spanning trees of
a graph, we make use of the well-known Matrix Tree Theorem 
(see e.g.~\cite{biggs}), as illustrated below.  
Given a directed graph $G$ with $n$ vertices, the negative Laplacian
of $G$ is the $n\times n$ matrix $L(G)$, where $L(G)_{v,w}$ ($v\neq
w$) equals the negative of the weight of the edge from vertex $v$ to
vertex $w$, and $L(G)_{v,v}$ is the weighted sum of the arcs emanating
from $v$.  The determinant of the submatrix obtained by deleting row
$r$ and column $r$ from $L(G)$ gives the (weighted) number of the
spanning trees rooted at $r$.

\centerline{\epsfig{figure=MTT.ps,width=\textwidth}}

We will make use of two ways to evaluate this determinant.  We can
exhibit $n-1$ orthogonal nonzero eigenvectors of the matrix obtained from
$L(G)$ by deleting row and column $r$, and multiply their eigenvalues.  
In the case where the graph is undirected,
so that the matrix is symmetric,
an alternate procedure is available:
we can exhibit $n$ orthogonal
nonzero eigenvectors of $L(G)$, 
multiply their eigenvalues except the
zero eigenvalue, and divide by $n$.  
(A more accurate general description of this procedure
is that the multiplicity of the eigenvalue zero
is to be reduced by 1.
When zero is a multiple eigenvalue,
the procedure gives rise to the product 0.
However, in all our examples, zero is a simple eigenvalue,
corresponding to the fact that the graph is connected
and hence possesses one or more spanning trees.)
This second method can be
viewed as a variation on the first method, 
in which an auxiliary vertex is added to the graph,
the edges from every other vertex to it are given weight $\varepsilon$,
the trees rooted at the auxiliary vertex are counted,
$\varepsilon$ is sent to 0,
and the result is divided by $n$,
corresponding to the fact that any of the $n$ vertices of the original graph
could become the unique vertex joined to the extra vertex
in a spanning tree in the new graph.
(The details are left to the reader.)

Every graph considered in the examples below
either is a finite induced subgraph of an infinite square grid
or else is obtained from such a graph by adding a single additional vertex.
In each case we give an explicit formulas for the eigenfunctions.
We hasten to say that for a general planar graph
an explicit diagonalization would be much more difficult.
The tractability of our chosen examples arises from their connection
with the negative Laplacian of the infinite square grid $\Z^2$.
A function $f$ on $\Z^2$ is an eigenfunction of the negative Laplacian
if it satisfies the equation
$$4 f(x,y) - f(x-1,y) - f(x+1,y) - f(x,y-1) - f(x,y+1) = \lambda f(x,y)$$
for all integers $x,y$,
in which case $\lambda$ is the associated eigenvalue.
For each pair of complex numbers $\zeta, \zeta'$ 
satisfying $|\zeta| = |\zeta'| = 1$,
we can construct such an eigenfunction by putting
$f(x,y) = \zeta^x {\zeta'}^y$
and
$\lambda = 4-\zeta-\zeta^{-1}-\zeta'-{\zeta'}^{-1}$.
In many cases, the restriction of such a function $f(x,y)$ 
to points $x,y$ lying in some finite region in $\Z^2$
is an eigenfunction of the matrix associated with that region
by the Matrix Tree Theorem.

As a preparatory example,
we briefly mention here the trivial case
of counting spanning trees of an undirected chain $C$
consisting of $n$ vertices and $n-1$ edges.
The space of complex-valued functions on the vertices of $C$
can be identified with the $n$-dimensional space $W$ 
of odd periodic functions of period $2n+2$,
i.e.\ functions $f: \Z \mapsto {\mathbb C}$
that satisfy $f(-x) = f(2n+2-x) = -f(x)$ for all $x$ in $\Z$
(and that consequently satisfy $f(0) = f(n+1) = \dots = 0$).
Under this identification, the negative Laplacian of $C$ is carried over
to the negative Laplacian of $\Z$, giving rise to an eigenbasis for $W$
of the form $f_k(x) = \sin \frac{kx\pi}{n+1}$ ($1 \leq k \leq n$).

We draw each graph (see Figures~\ref{fig:oddXodd-} through \ref{fig:explicit-hex})
so that the non-root vertices are located at points
of the lattice, and so that every edge of the graph connects nearest
neighbors in the lattice.  With high-degree root vertices it is
necessary to draw the root vertex in an extended fashion, covering
multiple points of the lattice, to ensure that all the edges incident
to the root are drawn between neighboring points in the lattice.

Generally the requisite eigenvectors of the matrix $L(G)$ (either with or
without the root row and column deleted) will be eigenvectors of the
negative Laplacian of the infinite lattice that also satisfy certain boundary
conditions (analogous to the conditions $f(0) = f(n+1) = 0$ in the
preparatory example).  Usually the infinite lattice will be $\Z^2$
although sometimes it will be $(\Z+\frac12)^2 =
 \{(x,y): x-\frac12,y-\frac12 \in \Z\}$.
There are two types of boundary conditions.  Suppose that $v$ and $w$
are nearest neighbors in the lattice, and $v$ is a non-root vertex of
the graph.  If $w$ is not in the graph, then we require that
$f(v)=f(w)$.  If the root vertex is drawn so as to contain the point
$w$, and the root row and column are deleted from $L(G)$,
then we require that $f(w)=0$.  We leave it to the reader to
check that if $f$ satisfies these boundary conditions, then the
restriction of $f$ to $G$ is an eigenvector of the matrix with
eigenvalue $\lambda$.

For certain subgraphs of the square-octagon graph,
exact formulas for the number of perfect matchings can be found by 
doing urban renewal to get a weighted version of one
of the graphs shown below.
In each case
the eigenvectors of the weighted version can be obtained from
the eigenvectors of the unweighted version  by multiplying by
weights $2^{i+j}$ and $2^{-i-j}$ 
as in \sref{square-octagon}.
The eigenvalue in
the weighted version is obtained from the unweighted eigenvalue by adding $1$.

In the rest of this section we give a number of graphs $G$, their duals
$G^\perp$, the graphs $\H(v^*,f^*)$, 
the eigenvectors of $L(G)$ (possibly with root row
and column removed), and the corresponding formula for the number of
perfect matchings of $\H(v^*,f^*)$.  When $G$ is undirected, we can also
enumerate the matchings of $\H(v^*,f^*)$ by counting the 
spanning trees of the dual of $G$.  

\subsection{Temperley's bijection}
Temperley's original bijection involved computing the number
of perfect matchings of a $2\ell-1$ by $2m-1$ subgraph of the square grid
with a corner removed.
By the main theorem, this is the number of spanning trees of an $\ell\times m$ 
rectangle.
\begin{figure}[ht]
\centerline{\epsfig{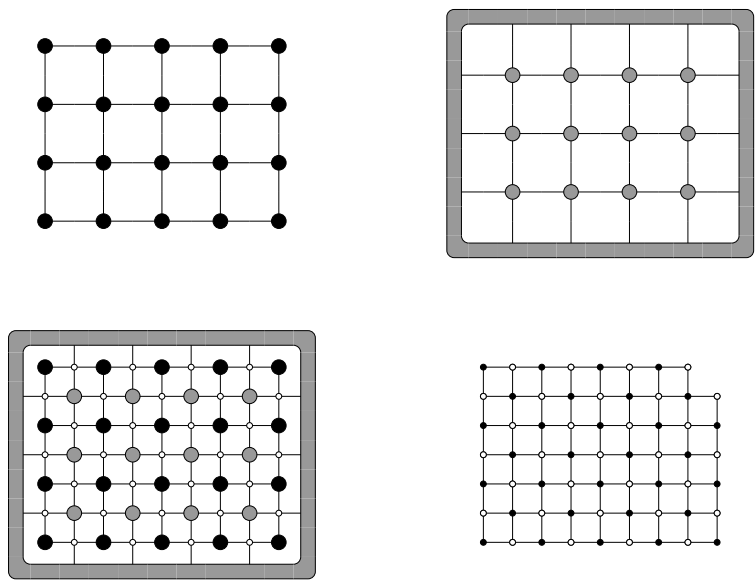}}
\caption{Temperley's original bijection, 
for the odd-by-odd rectangular region with a corner removed.
$\ell=5$, $m=4$.}
\label{fig:oddXodd-}
\end{figure}
The eigenvectors of the negative Laplacian are of the form
$$f(x,y)= f_{j,k}(x,y)=\cos\frac{\pi j x}{\ell}\cos\frac{\pi k y}{m},$$
where $x$ runs from $\frac12$ to $\ell-\frac12$ by integer steps
and $y$ from $\frac12$ to $m-\frac12$ by integer steps, 
$j$ is an integer in $[0,\ell-1]$ and $k$ is an integer in $[0,m-1]$.
(The lower left vertex is $(x,y)=(\frac12,\frac12)$ and the upper right
vertex is $(x,y)=(\ell-\frac12,m-\frac12$.)
The eigenvalue of $f_{j,k}$ is $4-2\cos\frac{\pi j}{\ell}-
2\cos\frac{\pi k}{m}$,
which is zero when $k=\ell=0$.
The number of spanning trees is
$$\frac1{\ell m}
\prod
\Bigl[4-2\cos\frac{\pi j}{\ell}-
2\cos\frac{\pi k}{m}\Bigr]$$
where the product is taken over all pairs $(j,k)$
with $0 \leq j \leq \ell-1$, $0 \leq k \leq m-1$
except $(0,0)$.

Using panel (b) of Figure~\ref{fig:oddXodd-},
we can compute the number of spanning trees
in a different way. In this case we take the negative Laplacian 
of the graph in panel (b) and remove a
row and column corresponding to the outer vertex. 
The eigenvectors of the resulting graph are as follows.
$$f_{j,k}(x,y)=\sin\frac{\pi j x}{\ell}\sin\frac{\pi k y}{m},$$
where this time $x$ runs from $0$ to $\ell$ by integer steps
and $y$ from $0$ to $m$ by integer steps, 
$j$ is an integer in $[1,\ell-1]$ and $k$ an integer in $[1,m-1]$.
Note that the $x,y$ coordinates correspond to the centers of the faces
in the coordinates of panel (a). 
The corresponding eigenvalue is $4-2\cos\frac{\pi j}{\ell}-
2\cos\frac{\pi k}{m}$.
The number of spanning trees is
$$
\prod_{j=1}^{\ell-1}\prod_{k=1}^{m-1}\Bigl[4-2\cos\frac{\pi j}{\ell}-
2\cos\frac{\pi k}{m}\Bigr].$$
In this case there is no zero eigenvalue since we have removed a row
and column of the negative Laplacian.
The equivalence of these two formulas follows from the 
well-known identity for $\ell=\prod_{j=1}^{\ell-1}(2-2\cos\frac{\pi j}{\ell})$.

\subsection{Dimers on an even-by-odd rectangle}
\begin{figure}[h]
\centerline{\epsfig{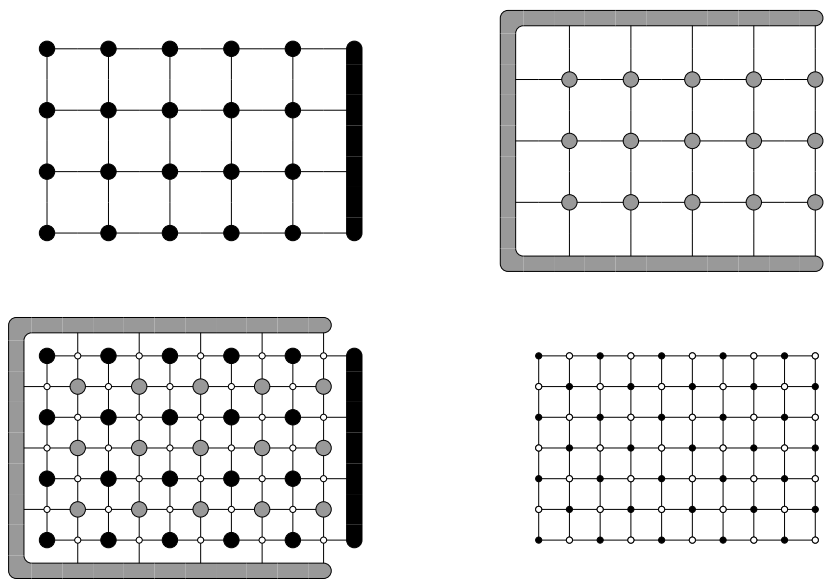}}
\caption{The even-by-odd rectangular region.  $\ell=5$, $m=4$.}
\label{fig:evenXodd}
\end{figure}
See \fref{evenXodd}.
In this case we again remove the row and column from the negative Laplacian
corresponding to the ``extended'' vertex.
The eigenvectors of the negative Laplacian are of the form
$$f_{j,k}(x,y)=\cos\frac{\pi(2j+1)x}{2\ell+1}\cos\frac{\pi k y}{m},$$
where $x$ runs from $\frac12$ to $\ell+\frac12$ and $y$ from 
$\frac12$ to $m-\frac12$, $j$ is an integer in $[0,\ell-1]$ and 
$k$ is an integer in $[0,m-1]$.
The corresponding eigenvalue is $4-2\cos\frac{\pi(2j+1)}{2\ell+1}-2
\cos\frac{\pi k}{m}$.
The number of spanning trees is
$$\prod_{j=0}^{\ell-1}\prod_{k=0}^{m-1}\Bigl[4-2\cos\frac{\pi (2j+1)}{2\ell
+1}-
2\cos\frac{\pi k}{m}\Bigr].$$

Similarly, in panel (b) with the extended vertex removed the eigenvectors
are
$$f_{j,k}(x,y)=\sin\frac{\pi(2j+1)x}{2\ell+1}\sin\frac{\pi k y}{m},$$
where $x$ runs from $0$ to $\ell$ and $y$ from 
$0$ to $m$, $j$ is an integer in $[0,\ell-1]$ and 
$k$ is an integer in $[1,m-1]$.
The corresponding eigenvalue is $4-2\cos\frac{\pi(2j+1)}{2\ell+1}-2
\cos\frac{\pi k}{m}$.

\subsection{Dimers on an even-by-even rectangle}
\begin{figure}[htbp]
\centerline{\epsfig{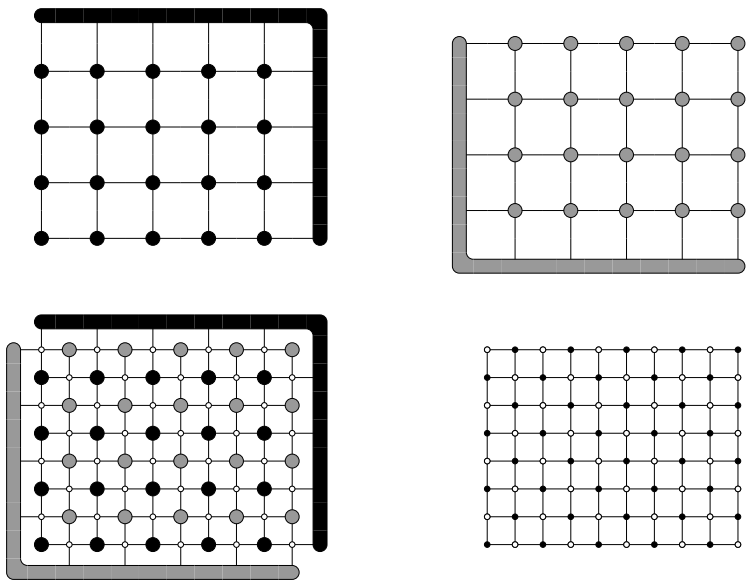}}
\caption{The even-by-even rectangular region.  $\ell=5$, $m=4$.}
\label{fig:evenXeven}
\end{figure}
See \fref{evenXeven}.
In panel (a) the eigenvectors of the negative Laplacian 
after removal of the extended vertex are
$$f_{j,k}(x,y)=\cos\frac{\pi(2j+1)x}{2\ell+1}\cos\frac{\pi(2k+1)y}{2m+1},$$
where $x$ runs from $\frac12$ to $\ell+\frac12$ and $y$ from 
$\frac12$ to $m+\frac12$, $j$ is an integer in $[0,\ell-1]$ and 
$k$ is an integer in $[0,m-1]$.
The number of spanning trees is
$$\prod_{j=0}^{\ell-1}\prod_{k=0}^{m-1}\Bigl[4-2\cos\frac{\pi (2j+1)}{2\ell
+1}-
2\cos\frac{\pi (2k+1)}{2m+1}\Bigr].$$

For eigenvectors in panel (b) replace the cosines with sines,
the $x$ range to $0$ to $\ell$, and the $y$ range $0$ to $m$.
The formula for the determinant is identical.

\subsection{Dimers on an odd-by-odd rectangle with an extra vertex}
\begin{figure}[tbp]
\centerline{\epsfig{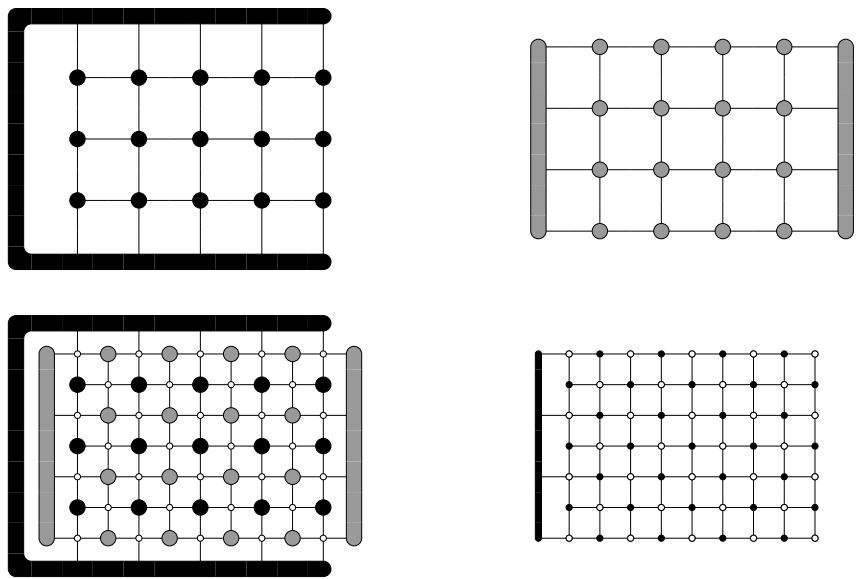}}
\caption{The odd-by-odd rectangular region with an extra vertex
connecting every other vertex on the left side.  $\ell=5$, $m=4$.}
\label{fig:oddXodd+}
\end{figure}

See \fref{oddXodd+}.
In panel (a) the eigenvectors of the negative Laplacian 
after removal of the extended vertex are
$$f_{j,k}(x,y)=\cos\frac{\pi jx}{\ell}\sin\frac{\pi ky}{m},$$
where $x$ runs from $\frac12$ to $\ell-\frac12$ and $y$ from
$0$ to $m$, $j$ is an integer in $[0,\ell-1]$ and 
$k$ is an integer in $[1,m-1]$.
The number of spanning trees is
$$\prod_{j=0}^{\ell-1}\prod_{k=1}^{m-1}\Bigl[4-2\cos\frac{\pi j}{\ell}-
2\cos\frac{\pi k}{m}\Bigr].$$

The eigenvectors in panel (b) are more complicated.

\smallskip

\begin{remark}
The formula for this region can also be derived by multiplying $m$
(the number of ways the extra vertex can be paired) by Temperley's
formula for the odd-by-odd region with a corner removed, since
Temperley's formula still holds with other perimeter vertices of the
right parity removed instead of the corner.
\end{remark}

\subsection{A diamond-shaped region}
\begin{figure}[tbp]
\centerline{\epsfig{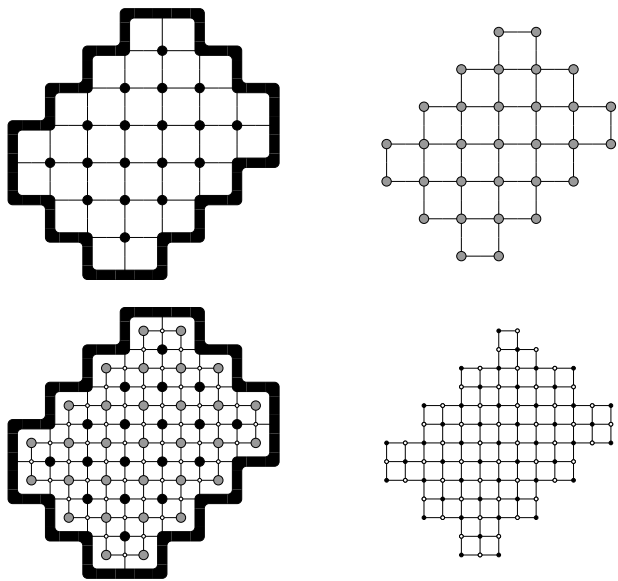}}
\caption{A diamond-shaped region with $\ell=4,m=3$.}
\label{fig:diamond}
\end{figure}
See \fref{diamond}.
In panel (a)
the eigenvectors of the negative Laplacian 
after removal of the extended vertex are
$$f_{j,k}(x,y)=\sin\frac{\pi j(x+y)}{2 \ell}\sin\frac{\pi k(x-y)}{2m},$$
where the origin $(x,y)=(0,0)$ is in the center bottom part of the extended
outer vertex of $G$, $x\in[-m,\ell]$, and $y\in[0,\ell+m]$, 
with $0\leq x+y\leq 2 \ell$ and $0\leq y-x \leq 2m$.
The indices $(j,k)\in[1,2 \ell-1]\times[1,m-1]\cup[1,\ell]\times\{m\}.$
The number of spanning trees is
$$\prod_{j,k}\Bigl[4-4\cos\frac{\pi j}{2 \ell}
\cos\frac{\pi k}{2m}\Bigr],$$
where $(j,k)$ runs over this range.

\subsection{Another diamond-shaped region}
\begin{figure}[htbp]
\centerline{\epsfig{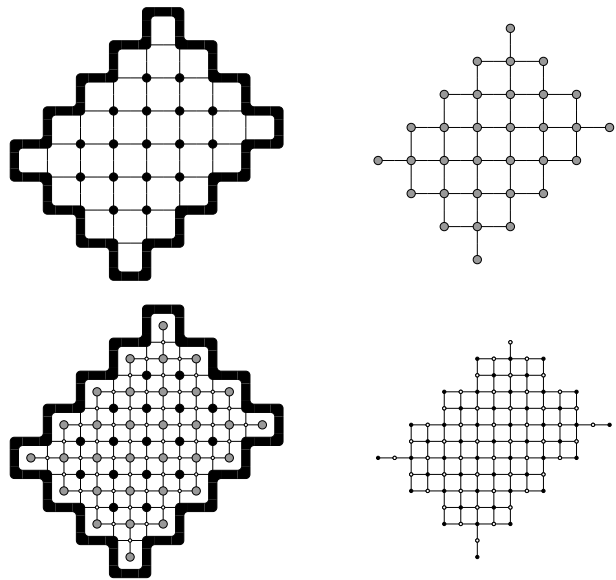}}
\caption{Another diamond-shaped region, with $\ell=4,m=3$.}
\label{fig:diamond2}
\end{figure}
See \fref{diamond2}.
In panel (a)
the eigenvectors of the negative Laplacian 
after removal of the extended vertex are
$$f_{j,k}(xy)=\sin\frac{\pi j(x+y)}{2 \ell}\sin\frac{\pi k(x-y)}{2m},$$
where the origin $(x,y)=(0,0)$ is at the bottom vertex of $G^\perp$,
$x\in[-m,\ell]$, and $y\in[0,\ell+m]$, with $0\leq x+y\leq 2 \ell$ and $0\leq y-x \leq 2m$.
The indices $(j,k)\in[1,2 \ell-1]\times[1,m-1]\cup[1,\ell-1]\times\{m\}.$
The number of spanning trees is
$$\prod_{j,k}\Bigl[4-4\cos\frac{\pi j}{2 \ell}
\cos\frac{\pi k}{2m}\Bigr],$$
where $(j,k)$ runs over this range. Note that this is exactly
the formula of the previous section, except for the index
$(j,k)=(\ell,m)$ whose eigenvalue is $4$. As a consequence the
number of spanning trees of this diamond is exactly one-fourth 
the number of spanning trees of the previous diamond.
This fact was first observed by \cite{stanley:problem} (when $\ell=m$),
and was first proved by \cite{knuth:trees}
and \cite{ciucu:cellular} (without assuming $\ell=m$);
a generalization was stated and proved by \cite{chow:trees}.

Again we cannot compute the eigenvectors in panel (b).

Similar formulas hold for two additional diamond-shaped regions with
boundaries intermediate between those shown in
Figures~\ref{fig:diamond} and \ref{fig:diamond2}.

\subsection{A triangular region}
\begin{figure}[htbp]
\centerline{\epsfig{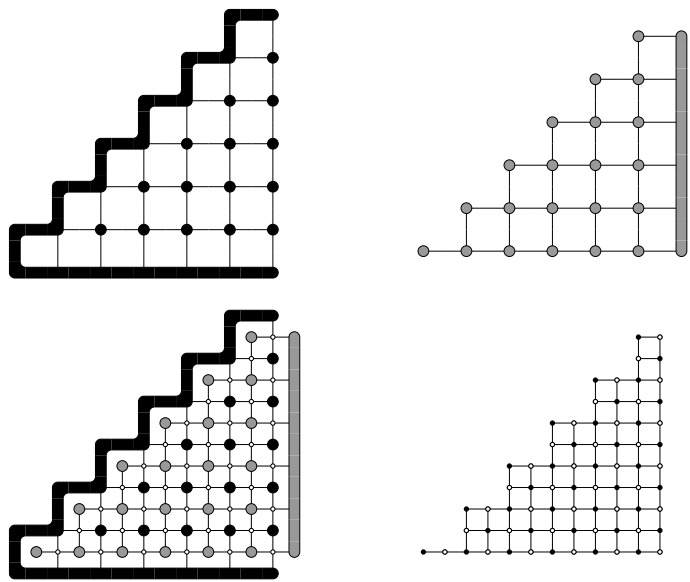}}
\caption{A triangular region.  $m=6$.}
\label{fig:triangle}
\end{figure}
See \fref{triangle}.
In panel (a) the eigenvectors of the negative Laplacian 
after removal of the extended vertex are
$$f_{j,k}(x,y)=\sin\frac{\pi (2j+1)x}{2m+1}\sin\frac{\pi (2k+1)y}{2m+1}
-\sin\frac{\pi (2k+1)x}{2m+1}\sin\frac{\pi (2j+1)y}{2m+1},$$
where the origin $(x,y)=(0,0)$ is the lower-bottom-most part of the extended
external vertex, $x$ and $y$ run from $0$ to $m$, with $x\geq y$,
and $j$ and $k$ are integers in $[0,m-1]$ with $j<k$.
The number of spanning trees is
$$\prod_{0\leq j<k<n}\Bigl[4-2\cos\frac{\pi (2j+1)}{2m+1}-2
\cos\frac{\pi (2k+1)}{2m+1}\Bigr].$$

We don't know how to compute the
eigenvectors of the negative Laplacian of the graph in panel (b).

\smallskip

\begin{remark}
This family of regions has been studied by Mihai Ciucu and Lior
Pachter.  \cite{ciucu:square} showed combinatorially that the number 
of domino tilings of the $2n\times 2n$ square equals $2^n$ times 
the square of the number of domino tilings of a particular region $H_n$, 
so one way to prove the claim (not the one given here!) is
to hitch a ride on the formula for the number of domino tilings of the
square --- or you could turn this around and derive the formula for the
number of tilings of the square as a corollary of the above formula.
Also, \cite{pachter:parity} gave a purely combinatorial 
proof that the number of domino tilings of $H_n$ is odd.
A different way to see this uses
2-adic analysis of the factors in the double product \citep{cohn:adic}.
\end{remark}

\subsection{Another triangular region, and the quartered Aztec diamond}
\begin{figure}[htbp]
\centerline{\epsfig{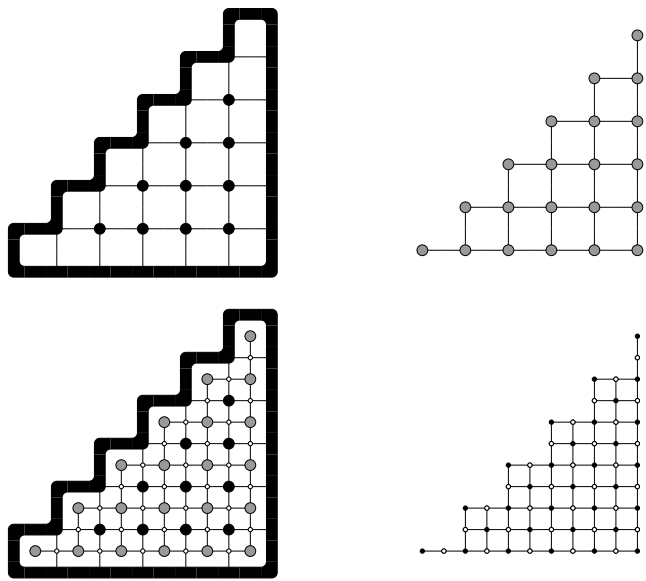}}
\caption{Another triangular region.  $G^{\perp}$ in panel (b) is the ``quartered Aztec diamond.''  $m=6$.}
\label{fig:triangle3}
\end{figure}
See \fref{triangle3}.

In panel (a) the eigenvectors of the negative Laplacian 
after removal of the extended vertex are
$$f_{j,k}(x,y)=\sin\frac{\pi j x}{m}\sin\frac{\pi k y}{m}
-\sin\frac{\pi k x}{m}\sin\frac{\pi j y}{m},$$
where the origin $(x,y)=(0,0)$ is the lower-bottom-most part of the extended
external vertex, $x$ and $y$ run from $0$ to $m$, with $x\geq y$,
and $j$ and $k$ are integers in $[1,m-1]$ with $j<k$.
The number of spanning trees is
$$\prod_{0<j<k<m}\Bigl[4-2\cos\frac{\pi j}{m}-2\cos\frac{\pi k}{m}\Bigr].$$

We don't know how to compute the
eigenvectors of the negative Laplacian of the graph in panel (b).

\smallskip

\begin{remark}
A formula for the number of spanning trees of the ``quartered Aztec
diamond'' ($G^\perp$ here) had been an open problem for some years
before Mihai Ciucu found an expression for it
(in work that has not yet been written up).
Here we effectively
get a formula for it by doing Fourier analysis on the dual graph.
The equivalence of Ciucu's formula and ours does not appear to be
an effortless identity.
\end{remark}

\subsection{A subregion of the hexagonal lattice}
\label{ssec:expl-hex}

\begin{figure}[htbp]
\centerline{\epsfig{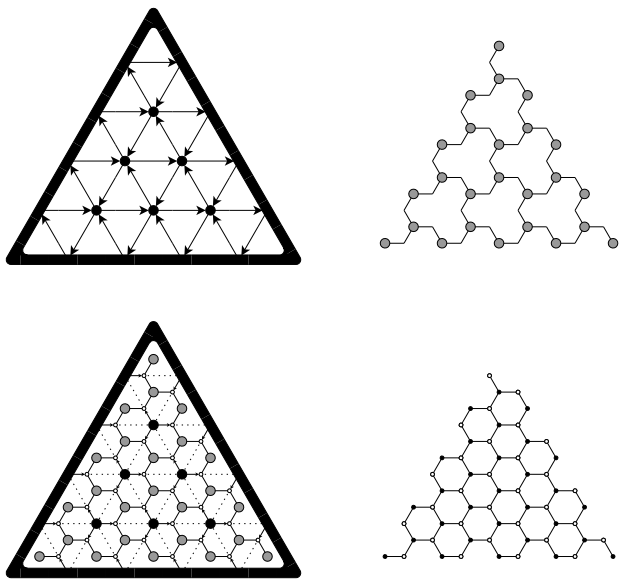}}
\caption{A triangular region of the hexagonal lattice.  Here $n=5$, and $m=n-2=3$.}
\label{fig:explicit-hex}
\end{figure}

Let $\omega=e^{2\pi i/3}$ and $\Z[\omega]$ be the lattice of Eisenstein
integers.  Let $G$ be the directed graph whose vertices are
$\Z[\omega]$, each vertex $v$ having an edge directed towards the
three vertices $v+1,v+\omega,v+\omega^{-1}$, as in
Figures~\ref{fig:hextree}a and~\ref{fig:explicit-hex}a.

Here we obtain an exact formula for the number of directed spanning
trees on the triangular region $T_m$ in $G$ whose vertices are
$0,m,-m\omega^{-1}$, with ``wired'' boundary conditions (\fref{explicit-hex}).

Let $n=m+2$.  Let $\Lambda_n$ be the sublattice of $\Z[\omega]$
generated by $n(1-\omega^{-1})$ and $n(\omega-\omega^{-1})$.  A
fundamental domain for $\Lambda_n$ is given by the hexagon whose
vertices are
$0,n,n(1-\omega^{-1}),-2n\omega^{-1},n(\omega-\omega^{-1}),n\omega$.
This fundamental domain is tiled with six copies of the equilateral
triangle whose vertices are $0,n,-n\omega^{-1}$.  The graph $T_m$ is
embedded in this triangle as the set of vertices not touching the
boundary.

Let $f_{\alpha,\beta,\gamma}$ be the function on $\Z[\omega]$ such
that $f_{\alpha,\beta,\gamma}(x+y\omega+z\omega^{-1}) = \alpha^x
\beta^y \gamma^z$; single-valuedness implies $\alpha \beta \gamma =
1$.  The function $f_{\alpha,\beta,\gamma}$ is an eigenvector of the
negative Laplacian of $G$, and has eigenvalue $3-\alpha-\beta-\gamma$.  In
order for $f_{\alpha,\beta,\gamma}$ to be a function on the torus
$G/\Lambda_n$, the relations $(\alpha/\gamma)^n=1$ and
$(\beta/\gamma)^n=1$ must also hold.  If $\alpha$ is a $3n$th root of
unity, and $\alpha/\gamma$ is an $n$th root of unity, then
$\beta/\gamma = 1/(\alpha\gamma^2)$ is also an $n$th root of unity, so
that $f_{\alpha,\beta,\gamma}$ is an eigenvector of the torus
$G/\Lambda_n$ with eigenvalue $3-\alpha-\beta-\gamma$.  The various
$f_{\alpha,\beta,\gamma}$'s formed in this way are orthogonal on the
torus, and since there are $3 n^2$ of them, the same as the number of
vertices in the torus, they form an eigenbasis of the negative Laplacian of the
torus.

Let $L_1,L_2,L_3$ be the three lines of symmetry of the lattice
$\Z[w]$ defined by $L_1=\R$, $L_2=\omega\R$, and $L_3=\omega^{-1}\R$.
Orthogonal reflection in each of these lines is not only a symmetry of
the lattice but preserves $\Lambda_n$ {\it and\/} the edge directions.
Therefore these reflections are symmetries of the underlying graph on
the torus.  The lines $L_1,L_2,L_3$ project to the torus, cutting it
into $6$ equilateral triangles.  Say that a function on the torus is
skew symmetric if reflecting it through any of the lines $L_1,L_2,L_3$
is equivalent to negating the function.  Since reflection in a line of
symmetry $L_1$, $L_2$, or $L_3$ sends the eigenvector
$f_{\alpha,\beta,\gamma}$ on the torus to the eigenvector
$f_{\alpha,\gamma,\beta}$, $f_{\gamma,\beta,\alpha}$, or
$f_{\beta,\alpha,\gamma}$ respectively (in each case two subscripts
have been transposed), when we express any skew symmetric function as
a linear combination of the $f_{\alpha,\gamma,\beta}$'s, we find that
it is a linear combination of the $g_{\alpha,\gamma,\beta}$'s, which
are defined by $$g_{\alpha,\beta,\gamma} = \sum_{\sigma\in S_3}
(-1)^\sigma f_{\sigma(\alpha,\beta,\gamma)} = f_{\alpha,\beta,\gamma}
- f_{\alpha,\gamma,\beta} + f_{\gamma,\alpha,\beta} -
f_{\gamma,\beta,\alpha} + f_{\beta,\gamma,\alpha} -
f_{\beta,\alpha,\gamma}.$$ $g_{\alpha,\beta,\gamma}$ is an eigenvector
of the torus $G/\Lambda_n$ with eigenvalue $3-\alpha-\beta-\gamma$.
If we restrict our attention to triples $(\alpha,\beta,\gamma)$
satisfying $\alpha<\beta<\gamma$ under some arbitrary ordering of the
roots of unity, then the $g_{\alpha,\beta,\gamma}$'s still span the
space of skew symmetric functions; furthermore a nontrivial linear
relation among them would imply a nontrivial linear relation among the
$f_{\alpha,\beta,\gamma}$'s, so the $g_{\alpha,\beta,\gamma}$'s such
that $\alpha<\beta<\gamma$ form an eigenbasis of the space of skew
symmetric functions.

Since the $g_{\alpha,\beta,\gamma}$'s are zero on the lines of
symmetry $L_1,L_2,L_3$, when we restrict the
$g_{\alpha,\beta,\gamma}$'s to one of the $6$ equilateral triangles,
they form an eigenbasis of the matrix obtained from the negative Laplacian of
the triangular subgraph by deleting the vertices on the boundary.
Thus multiplying their eigenvalues gives the number of spanning trees
of the triangular region with wired boundary conditions, rooted at the
wired boundary:
$$ \prod_{\substack{
          \alpha,\beta,\gamma \\
          \alpha\beta\gamma=1 \\
          \alpha^{3n} = 1     \\
          (\alpha/\beta)^n=1  \\
          \text{$\alpha,\beta,\gamma$ distinct}}}
       (3-\alpha-\beta-\gamma)^{1/6}.$$

The prime factors of the numbers given by this formula exhibit
some rather nice patterns.

\section{Open Problems} \label{sec:open}

Our work suggests (or might be useful in the solution of) 
a number of different problems.

First, there is the natural question of extending the correspondence
between spanning trees and matchings so that it applies to matchings
of more general planar graphs (or perhaps even some non-planar ones).
An intriguing example is the ``12-6-4 lattice'' (the infinite graph
obtained by taking the 1-skeleton of the Archimedean tiling of the
plane by dodecagons, hexagons, and squares) shown in \fref{12-6-4}.
\begin{figure}[htb]
\centerline{\epsfig{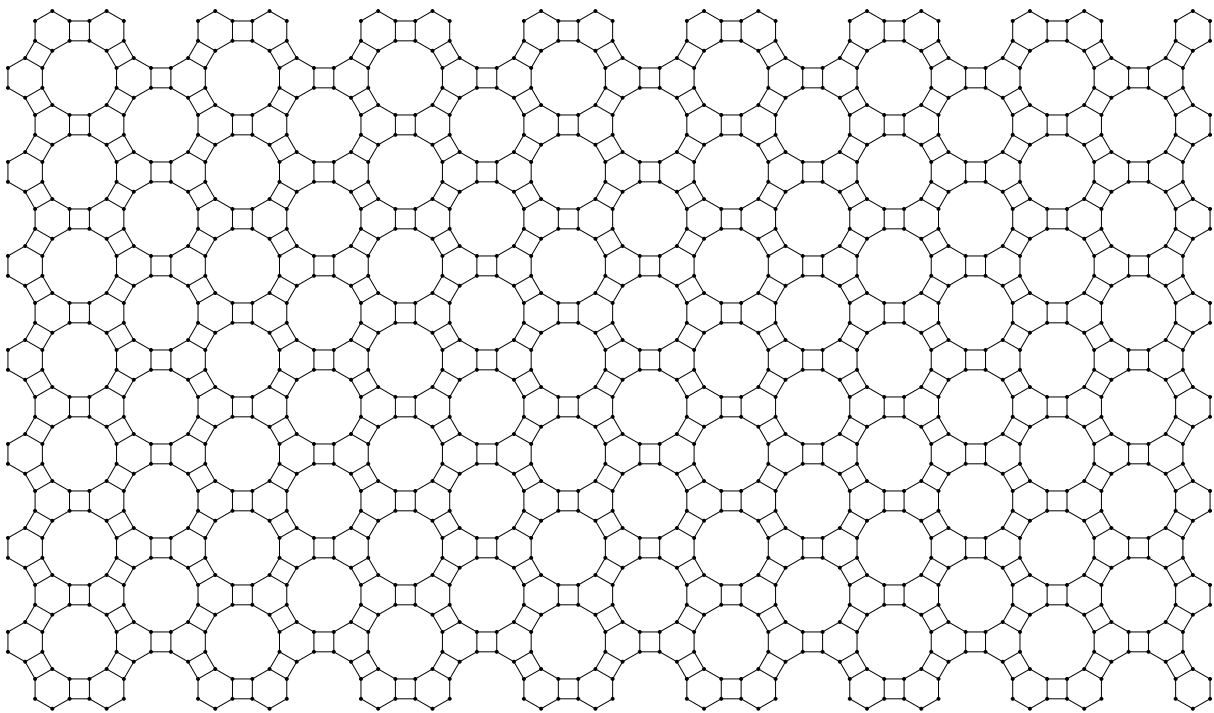}}
\caption{A portion of the ``12-6-4 lattice'' of dodecagons, hexagons, and
 squares.}
\label{fig:12-6-4}
\end{figure}
In a matching, each
vertex is paired to another vertex via one of three types of edges:
a 6/4 edge (bordering a hexagon
and square), a 12/4 edge (bordering a dodecagon and square), or a
12/6 edge (bordering a dodecagon and hexagon).  
In random perfect matchings of suitably defined
subgraphs (chosen so as to minimize the effect of the boundary),
the probabilities of
these three events are respectively $1/6+(19/78)R$, $1/3+(1/39)R$, and
$1/2-(7/26)R$, where $R$ is given by
$$ R = \sum_{n=0}^{\infty} \sum_{k=0}^n \binom{2 n}{n} \binom{n}{k}^2 \frac{14^k}{13^{2 n}}.$$
Curiously, these probabilities also turn up in the random spanning
trees of a certain directed weighted cartesian lattice.  In this
lattice each leftward edge has weight 25, each upward edge has weight
25, each rightward edge has weight 14, and each downward edge has
weight 1.  The probability that in the tree the parent of a vertex is
to the left is the same as the probability of a 12/4 edge in the
12-6-4 lattice, and the probability that the parent is above the given
vertex is the same as the probability of a 6/4 edge in the 12-6-4
lattice.  These ``coincidences'' suggest a weighted bijective
correspondence analogous to the one for the lattice of octagons and
squares, but we have been unable to find one.

One side-issue of a number-theoretic nature concerns the formula
given at the end of \S~\ref{ssec:expl-hex}.  As was mentioned, there are 
some interesting patterns governing the prime factors of these
numbers; for instance, it appears that most of the large prime
factors, including the very largest, are congruent to 1 mod $n$. 
It would be good to know why this is true.

Turning to entropy,
it would also be desirable to have a more general understanding
of asymptotics.
In all the examples considered in \sref{formula},
the logarithm of the number of spanning trees,
divided by the size of the system,
tends to a single limit that is independent
of the shape of the boundary,
but is a numerical characteristic of the infinite square grid.
One might ask the same question for more general grids.
For example, suppose we have some locally finite infinite planar graph $G$
that admits an adjacency-preserving action of $\Z^2$
that has only finitely many vertex-orbits and edge-orbits.
What is the limit of the normalized logarithm of
the number of spanning trees of large finite subgraphs of $G$,
as the subgraphs grow to fill $G$?
It is not too hard to guess the shape of the answer,
using the spectrum of the negative Laplacian,
but justifying this answer would require some care.

There is also the issue of uniqueness of the measure of maximum entropy.
In the case of the square grid, \cite{burton-pemantle:tree}
proved that there is a unique translation-invariant measure on spanning
trees of the infinite square grid that achieves the entropy bound.
This yields an analogous result on domino tilings of the plane.  It is
natural to ask the question for lozenge tilings of the plane; by the
theorems of this article, this is closely related to the problem of
determining whether there is a unique measure of maximum entropy for
directed spanning trees in the directed triangular lattice.

The Ising partition function and the number of spanning trees of a
graph are both evaluations of the Tutte polynomial of the graph at
special points; see \cite{welsh:tutte} for background on the Tutte
polynomial.  Evaluating the Tutte polynomial at most points is
\#P-complete, but at special points the polynomial can be evaluated in
polynomial time.  For planar graphs these special points include the
values for the Ising partition function and the number of spanning
trees.  For planar graphs there are bijective connections between
Ising systems and perfect matchings \citep{fisher:ising-dimer}, and
between spanning trees and perfect matchings.  Thus it is natural to
search for connections between perfect matchings and the other special
points of the Tutte polynomial.

\section*{Acknowledgements}
Part of this research was done while the first author was visiting Microsoft Research.
The second author thanks Noam Elkies, Sergey Fomin, and Greg Kuperberg
for useful conversations on the topic.

\bibliography{eus,kpw}
\bibliographystyle{plainnat}

\end{document}